\documentclass[11pt,a4paper]{article}

\usepackage{theorem,enumerate}
\usepackage{amsmath,latexsym,amssymb,amsfonts}
\usepackage{mathrsfs}

\newcounter{Scounter}
\def\sco{\arabic{Scounter}}
\def\ssco{\stepcounter{Scounter}}
\def\sets{\setcounter{Scounter}{1}}

\theorembodyfont{\normalfont\slshape}
\newtheorem{thm}{Theorem}

\newtheorem{Thm}{Theorem}
\renewcommand{\theThm}{\Alph{Thm}}

\newtheorem{lem}{Lemma}

\newtheorem{Q}{Question}
\renewcommand{\theQ}{}
\newtheorem{claim}{Claim}

\newtheorem{con}[thm]{Conjecture}

\renewcommand{\theCon}{\Alph{Con}}
\newtheorem{Lem}{Lemma}
\renewcommand{\theLem}{\Alph{Lem}}
\newtheorem{prob}{\normalfont\itshape Problem}
\newtheorem{remark}{\normalfont\itshape Remark}

\newcommand{\proof}{\medbreak\noindent\textit{Proof.}\quad}
\newcommand{\Proof}{\medbreak\noindent\textbf{Proof.}\quad}
\newcommand{\PROOF}[1]{\medbreak\noindent\textbf{Proof #1.}\ }
\newcommand{\proofA}{\medbreak\noindent\textit{Proof of Theorem A.}\quad}

\newcommand{\case}[1]{\noin{\it Case #1}}
\newcommand{\Case}[1]{\noin{\bf Case #1}}
\newcommand{\sscase}[2]{\noin{\bf (#1) }#2}
\newcommand{\qed}{{$\quad\square$\vs{3.6}}}
\newcommand{\Qed}{{$\;\square$\vs{3.6}}}
\newcommand{\qqed}{{$\quad\square \square$\vs{3.6}}}

\newcommand{\en}{\text{\upshape end}}
\newcommand{\diff}{\text{\upshape diff}}
\newcommand{\dist}{\text{\upshape dist}}
\newcommand{\diam}{\text{\upshape diam}}
\newcommand{\cl}{\text{\upshape cl}}
\newcommand{\bss}{\setminus}

\newcommand{\uppercut}[1]{\left\lceil #1 \right\rceil}
\newcommand{\vs}[1]{\vspace*{#1 mm}}
\newcommand{\hs}[1]{\hspace*{#1 mm}}
\newcommand{\noin}{\noindent}
\newcommand{\bs}{\setminus}
\newcommand{\kuu}[1]{{#1}'}
\newcommand{\hasi}[1]{\tilde{#1}}
\newcommand{\hasis}[1]{\tilde{#1}}
\newcommand{\dg}[2]{d_{#1}(#2)}
\newcommand{\bun}[2]{\frac{#1}{#2}}
\newcommand{\IRO}[2]{#1\ora{C}#2^-}
\newcommand{\OIR}[2]{#1^+\ora{C}#2}
\newcommand{\OIRO}[2]{#1^+\ora{C}#2^-}
\newcommand{\PR}[2]{P[#1,#2]}
\newcommand{\IR}[2]{#1\ora{C}#2}
\newcommand{\ep}{\varepsilon}
\newcommand{\ora}{\overrightarrow}
\newcommand{\ola}{\overleftarrow}
\newcommand{\el}{l}

\numberwithin{equation}{section}

\def\Vec#1{\mbox{\boldmath $#1$}}
\def\A{{ \mathcal{A}}}
\def\T{{ \mathcal{T}}}
\def\P{{ \mathcal{P}}}
\def\F{{ \mathcal{F}}}
\def\calH{{ \mathcal{H}}}
\def\C{{ \mathcal{C}}}
\def\I{{ \mathcal{I}}}
\def\Q{{ \mathcal{Q}}}
\def\T{{ \mathcal{T}}}
\def\G{{ \mathcal{G}}}
\def\X{{ \mathcal{X}}}
\def\U{{ \mathcal{U}}}
\def\R{{ \mathcal{R}}}

\makeatletter
\def\thanks#1{%
   \footnotemark
   \edef\@tempa{\noexpand\noexpand\noexpand\footnotetext[\the\c@footnote]}%
   \toks@\expandafter{\@thanks}%
   \toks\tw@{{#1}}
   \xdef\@thanks{\the\toks@\@tempa\the\toks\tw@}}
\makeatother

\def\ia{{\accent 19 \char 16}~}

\renewcommand{\baselinestretch}{1.3}

\bfseries\normalfont

\begin{document}

\title{Forbidden pairs and the existence of a dominating cycle}

\author{
Shuya Chiba$^{1}$\thanks{This work was partially supported by JSPS KAKENHI grant 23740087\\
\hspace{+14pt}
E-mail address: \texttt{schiba@kumamoto-u.ac.jp}}
\and
Michitaka Furuya$^{2}$\thanks{E-mail address: \texttt{michitaka.furuya@gmail.com} } 
\and 
Shoichi Tsuchiya$^{2}$\footnote{E-mail address: \texttt{wco.liew6.te@gmail.com}}
\vspace{+8pt}
 \\
\small
$^1$\small\textsl{Department of Mathematics and Engineering, 
Kumamoto University,}\\ 
\small\textsl{2-39-1 Kurokami, Kumamoto 860-8555, Japan}\\
\small
$^{2}$\small\textsl{Department of Mathematical Information Science,
Tokyo University of Science,}\\ 
\small\textsl{1-3 Kagurazaka, Shinjuku-ku, Tokyo 162-8601, Japan}\\
}

\date{}
\maketitle

\vspace{-20pt}
\begin{abstract}
A cycle $C$ in a graph $G$ is called dominating 
if every edge of $G$ is incident with a vertex of $C$. 
For a set $\calH$ of connected graphs, 
a graph $G$ is said to be $\calH$-free 
if $G$ does not contain any member of $\calH$ as an induced subgraph. 
When $|\calH| = 2$, $\calH$ is called a forbidden pair. 
In this paper, 
we investigate the set $\Vec{H}$ of pairs $\cal{H}$ of connected graphs which satisfies that 
every $2$-connected $\calH$-free graph has a dominating cycle.
In particular, we show that $\Vec{H}$ is a very small class of pairs of graphs 
and find some pairs of graphs which belong to $\Vec{H}$.

\end{abstract}

\section{Introduction}
\label{introduction}

A cycle $C$ in a graph $G$ is called \textit{dominating} if every edge of $G$ is incident with a vertex of $C$.
In this paper, we investigate forbidden subgraphs which imply the existence of a dominating cycle.
The origin of our research goes back to results on forbidden subgraphs implying the existence of a Hamilton cycle.

All graphs considered here are finite simple graphs. 
For standard graph-theoretic terminology not explained in this paper, we refer the readers to \cite{D}. 
A graph $G$ is said to be \textit{Hamiltonian} if $G$ has a \textit{Hamilton cycle}, i.e., a cycle containing all vertices of $G$.
The study on a Hamilton cycle is one of the most important and basic topics in graph theory.
It is known that the problem of determining whether a given graph is Hamiltonian or not 
belongs to the class of $\Vec{NP}$-complete problems, that is, a difficult problem in a combinatorial sense. 
So, many researchers have studied sufficient conditions for Hamiltonicity of graphs, 
and 
there is a large amount of literature concerning conditions in terms of 
order, size, vertex degrees, independence number, forbidden subgraphs and so on 
(see a survey \cite{Br}).

Let $\calH$ be a set of connected graphs. 
A graph $G$ is said to be \textit{$\calH$-free} 
if $G$ does not contain $H$ as an induced subgraph for all $H$ in $\calH$, 
and we call each graph $H$ of $\calH$ a \textit{forbidden subgraph}. 
If $\calH = \{H\}$, 
then we simply say that $G$ is $H$-free. 
We call $\calH$ a \textit{forbidden pair} if $|\calH| = 2$. 
When we consider $\calH$-free graphs, 
we assume that each member of $\calH$ has order at least $3$ 
because $K_{2}$ is the only connected graph of order $2$ 
and connected $K_{2}$-free graphs are only $K_{1}$ 
(here $K_{n}$ denotes the complete graph of order $n$).
In order to state results on forbidden subgraphs clearly, 
we further introduce several notations. 
For two graphs $H_{1}$ and $H_{2}$, 
we write $H_{1} \prec H_{2}$ if $H_{1}$ is an induced subgraph of $H_{2}$,  
and for two sets $\calH_{1}$ and $\calH_{2}$ of connected graphs, 
we write $\calH_{1} \le \calH_{2}$ 
if for every graph $H_{2}$ in $\calH_{2}$, 
there exists a graph $H_{1}$ in $\calH_{1}$ with $H_{1} \prec H_{2}$. 
By the definition of the relation ``$\leq $'', if $\calH_{1} \le \calH_{2}$, then $\calH_{1}$-free graph is also $\calH_{2}$-free.

\begin{figure}
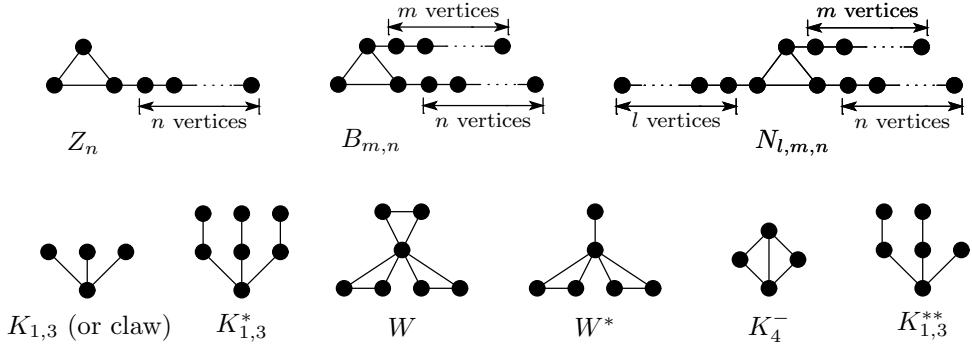

\begin{center}
\unitlength 0.1in
%
\caption{The forbidden subgraphs}
\label{forbidden subgraphs}
\end{center}
\end{figure}

The forbidden pairs that force the existence of a Hamilton cycle in $2$-connected graphs had been studied in \cite{BV,DGJ,GJ}.
Eventually, a characterization of such pairs was accomplished in \cite{Be} as follows
(here let $P_{n}$ denote the path of order $n$, 
and the graphs $K_{1, 3}$ (or claw), 
$B_{m, n}$ and $N_{l, m, n}$ 
are the ones that are depicted in Figure~\ref{forbidden subgraphs}).

\begin{Thm}[Bedrossian~\cite{Be}]
\label{Be}
Let $\calH$ be a set of two connected graphs. 
Then every $2$-connected $\calH$-free graph is Hamiltonian 
if and only if 
$\calH \le \{K_{1, 3}, P_{6}\}$,  $\calH \le \{K_{1, 3}, B_{1, 2}\}$, or $\calH \le \{K_{1, 3}, N_{1, 1, 1}\}$. 
\end{Thm}

On the other hand, Faudree, Gould, Ryj\'{a}\v{c}ek and Schiermeyer~\cite{FGRS} proved that
every $2$-connected $\{K_{1,3},Z_{3}\}$-free graph of order at least $10$ is Hamiltonian 
(here $Z_{n}$ is the one that is depicted in Figure~\ref{forbidden subgraphs}).
In \cite{FG}, 
the forbidden pairs for Hamiltonicity of $2$-connected graphs have been completely determined 
even when we allow a finite number of exceptions.

\begin{Thm}[Faudree and Gould~\cite{FG}]
\label{FG}
Let $\calH$ be a set of two connected graphs. 
Then every $2$-connected $\calH$-free graph of sufficiently large order is Hamiltonian 
if and only if 
$\calH \le \{K_{1, 3}, P_{6}\}$, 
$\calH \le \{K_{1, 3}, Z_{3}\}$, 
$\calH \le \{K_{1, 3}, B_{1, 2}\}$, 
or 
$\calH \le \{K_{1, 3}, N_{1, 1, 1}\}$. 
\end{Thm}

A \textit{$2$-factor} of a graph $G$ is a spanning subgraph of $G$ in which every component is a cycle.
It is known that a $2$-factor is one of the relaxed structures of a Hamilton cycle 
since a Hamilton cycle is a connected 2-factor. 
In fact, 
the sufficient conditions for the existence of a $2$-factor have been extensively studied 
in order to 
investigate the difference between the existence of a Hamilton cycle and a $2$-factor in graphs 
(see a survey \cite{G}).
As part of it,
the forbidden pairs that imply a $2$-connected graph has a $2$-factor 
was characterized by J.R.~Faudree, R.J.~Faudree and Ryj\'{a}\v{c}ek~\cite{FFR}.

\begin{Thm}[J.R.~Faudree, R.J.~Faudree and Ryj\'{a}\v{c}ek~\cite{FFR}]
\label{FFR}
Let $\calH$ be a set of two connected graphs. 
Then the following hold. 
\vspace{-5pt}
\begin{enumerate}[{\upshape(i)}]
\item Every $2$-connected $\calH$-free graph has a $2$-factor 
if and only if 
$\calH \le \{K_{1, 3}, P_{7}\}$, 
$\calH \le \{K_{1, 3}, Z_{4}\}$, 
$\calH \le \{K_{1, 3}, B_{1, 3}\}$, 
 or 
$\calH \le \{K_{1, 3}, N_{1, 1, 2}\}$. 
\vspace{-5pt}
\item Every $2$-connected $\calH$-free graph of sufficiently large order has a $2$-factor 
if and only if 
$\calH \le \{K_{1, 3}, P_{7}\}$, 
$\calH \le \{K_{1, 3}, Z_{4}\}$ 
$\calH \le \{K_{1, 3}, B_{1, 4}\}$, 
$\calH \le \{K_{1, 3}, N_{1, 1, 3}\}$,  
or 
$\calH \le \{K_{1, 4}, P_{4}\}$. 
\end{enumerate}
\end{Thm}

On the other hand, 
one often try to find a dominating cycle in order to find a Hamilton cycle in a given graph 
(recall that a cycle $C$ in a graph $G$ is dominating if every edge of $G$ is incident with a vertex of $C$). 
For example, if some longest cycle in a graph $G$ is dominating and the independence number of $G$ is at most its minimum degree, 
then $G$ has a Hamilton cycle (the related results can be found in \cite{BHJV, T}).
It is also shown that 
the dominating cycle conjecture that ``every cyclically $4$-edge-connected cubic graph has a dominating cycle'' by Fleischner~\cite{Fl}
is equivalent to not only the well-known conjecture that ``every $4$-connected $K_{1, 3}$-free graph is Hamiltonian'' by Matthews and Sumner~\cite{MS}
but also many other statements on Hamiltonicity of graphs 
(see a survey \cite{BRV}).
In this sense, a topic on a dominating cycle is one of important relaxations of a Hamilton cycle.

In this paper, our motivation is to investigate the difference between the existence of a Hamilton cycle and a dominating cycle 
of a $2$-connected graph in terms of the forbidden pair.

\begin{prob}
\label{determine H}
Determine the set $\Vec{H}$ (resp., $\Vec{H}'$) of pairs $\calH$ of connected graphs which satisfy that
every $2$-connected $\calH$-free graph 
(resp., every $2$-connected $\calH$-free graph of sufficiently large order) has a dominating cycle.
\end{prob}

Concerning the above problem,
we first show that $\Vec{H}$ and $\Vec{H}'$ are very small classes of pairs.
Let $K_{1, 3}^{*}$, $W$, $W^{*}$ and $K_{4}^{-}$ 
be the ones that are depicted 
in Figure~\ref{forbidden subgraphs}, 
and set 
$\calH_{1} = \{ K_{1, 3}, Z_{4} \}$, 
$\calH_{2} = \{ K_{1, 3}, B_{1, 2} \}$, 
$\calH_{3} = \{ K_{1, 3}, N_{1, 1, 1} \}$, 
$\calH_{4} = \{ P_{4}, W \}$, 
$\calH_{5} = \{ K_{1, 3}^{*}, Z_{1} \}$, 
$\calH_{6} = \{ P_{5}, W^{*} \}
$ 
and 
$
\calH_{7} = \{ P_{5}, K_{4}^{-} \}$.

\begin{thm}
\label{2-conn (necessity)}
Let $\cal{H}$ be a set of two connected graphs. 
If 
there exists a positive integer $n_{0} = n_{0}(\calH)$ such that 
every $2$-connected $\cal{H}$-free graph of order at least $n_{0}$ has a dominating cycle, 
then 
$\calH \le \calH_{i}$ for some $i$ with $1 \le i \le 7$. 
\end{thm}

Theorem~\ref{2-conn (necessity)} implies that 
$\Vec{H} \subseteq \Vec{H}' \subseteq 
\{ \calH : |\calH| = 2, \calH \le \calH_{i} \textup{ for some } i \textup{ with } 1 \le i \le 7 \} =: \Vec{H}^{*}$ 
So, 
the remaining problem is 
that 
whether 
$\calH \in \Vec{H}$ or not 
($\calH \in \Vec{H}'$ or not)
when $\calH$ is a member of $\Vec{H}^{*}$. 
We actually guess that 
the contrary of Theorem~\ref{2-conn (necessity)} holds.

\begin{con}
\label{con}
Let $\cal{H}$ be a set of two connected graphs. 
If $\calH \in \{\calH_{i} : 1 \le i \le 7\}$, 
then 
every $2$-connected $\calH$-free graph (of sufficiently large order) has a dominating cycle. 
\end{con}

As a partial solution of Conjecture~\ref{con}, 
in this article, 
we further show that 
$\calH_{i}$ is a member of $\Vec{H} \ ( \ \subseteq \Vec{H}')$ for $1 \le i \le 4$ 
and 
that a pair $\calH$ of connected graphs with 
$\calH \le \calH_{5}$ and $\calH \neq \calH_{5}$ 
is also a member of $\Vec{H} \ ( \ \subseteq \Vec{H}')$
(here $K_{1, 3}^{**}$ is the graph obtained from $K_{1, 3}^{*}$ by deleting one leaf
(see Figure~\ref{forbidden subgraphs})
and $\calH_{5}' = \{ K_{1, 3}^{**}, Z_{1} \}$).

\begin{thm}
\label{2-conn (sufficiency)}
Let $\cal{H}$ be a set of two connected graphs. 
If $\calH \in \big\{ \calH_{i} : 1 \le i \le 4 \big\} \cup \{ \calH_{5}' \}$, 
then 
every $2$-connected $\cal{H}$-free graph has a dominating cycle. 
\end{thm}

We prove Theorem~\ref{2-conn (necessity)} in Section~\ref{proof of necessity}, 
and slightly stronger statements than Theorem~\ref{2-conn (sufficiency)} in Section~\ref{proof of sufficiency}
(see Theorems~\ref{K_{1,3}-free}--\ref{K_{1,3}^{**}, Z_{1}}).

\begin{remark}
By observing Theorems~\ref{Be}, \ref{FG} and \ref{FFR}, 
one may think that we always need an induced subgraph of a star in forbidden pairs 
for Hamiltonicity-like properties of graphs.
In fact, as one of the approach to attack Matthews-Sumner conjecture,
the forbidden pair containing $K_{1, 3}$ 
for the existence of a Hamilton cycle in $k$-connected graphs ($k \ge 3$) have been also studied, e.g., see \cite{Fu, LXYY, LP, P}.
However, when we consider the existence of a dominating cycle, 
the situation is a bit different from Theorems~\ref{Be}, \ref{FG} and \ref{FFR},
i.e., there exist forbidden pairs which contain no star and force the existence of a dominating cycle in $2$-connected graphs (see Theorem~\ref{2-conn (sufficiency)}).
\end{remark}

\section{Terminology and notation}
\label{terminology}

In this section, 
we prepare terminology and  notation which we use in subsequent sections. 

Let $G$ be a graph. 
We denote by $V(G)$, $E(G)$ 
and $\Delta(G)$ 
the vertex set, the edge set 
and the maximum degree of $G$, respectively. 
For $X \subseteq V(G)$, 
we let $G[X]$ denote the subgraph induced by $X$ in $G$, 
and let $G - X = G[V(G) \bss X]$. 
Let $v$ be a vertex of $G$. 
We denote by $N_{G}(v)$ and $d_{G}(v)$ the neighborhood 
and the degree of $v$ in $G$, respectively.
For $X \subseteq V(G) \bss \{v\}$, 
we let $N_{G}(v; X) = N_{G}(v) \cap X$, 
and for $V, X \subseteq V(G)$ with $V \cap X = \emptyset$, 
let $N_{G}(V; X) = \bigcup_{v \in V}N_{G}(v; X)$.
We often identify a subgraph $F$ of $G$ with its vertex set $V(F)$
(for example, $N_{G}(v; V(F))$ is often denoted by $N_{G}(v; F)$). 
For a positive integer $l$, 
we define $V_{l}(G) = \{v \in V(G) : d_{G}(v) = l\}$. 
For $u, v \in V(G)$, 
$\dist_{G}(u, v)$ denotes the distance between $u$ and $v$ in $G$, 
and we define the \textit{diameter} $\diam(G)$ of $G$ by 
$\diam(G) = \max \{ \dist_{G}(u, v) : u, v \in V(G) \}$. 
When $G$ has a cycle, 
we denote by $c(G)$ the \textit{circumference} of $G$, i.e., 
the length of the longest cycle of $G$. 
A path with end vertices $u$ and $v$ 
is denoted by a \textit{$(u, v)$-path}. 

We write a cycle (or a path) $C$ with a given orientation by $\ora{C}$.
If there exists no chance of confusion, we abbreviate $\ora{C}$ by $C$. 
Let $\ora{C}$ be an oriented cycle or a path.
For $x,y \in V(C)$,
we denote by $x\ora{C}y$
the $(x, y)$-path on $\ora{C}$.
The reverse sequence
of $x\ora{C}y$
is denoted by $y\ola{C}x$.
For $u \in V(C)$, 
we denote the $h$-th successor and the $h$-th predecessor of $u$ 
on $\ora{C}$ by $u^{+h}$ and $u^{-h}$, respectively, 
and let $u^{+0} = u$. 
For $X \subseteq V(C)$, 
we define $X^{+h} = \{x^{+h} : x \in X\}$ 
and $X^{-h} = \{x^{-h} : x \in X\}$, respectively.  
We abbreviate $u^{+1}$, $u^{-1}$, $X^{+1}$ and $X^{-1}$ by $u^{+}$, $u^{-}$, $X^+$ and $X^-$, respectively.

For two graphs $G_{1}$ and $G_{2}$ with $V(G_{1}) \cap V(G_{2}) = \emptyset$, 
let $G_{1} \cup G_{2}$ denote the union of $G_{1}$ and $G_{2}$, 
and let $G_{1} + G_{2}$ denote the join of $G_{1}$ 
and $G_{2}$, i.e., the graph obtained from $G_{1} \cup G_{2}$ 
by joining each vertex in $V(G_{1})$ to all vertices in $V(G_{2})$. 
For a graph $G$ and $l \ge 1$, let $lG$ denote the union of 
$l$ vertex-disjoint copies of $G$.

\section{Proof of Theorem~\ref{2-conn (necessity)}}
\label{proof of necessity}

We first prepare the following lemma concerning 
the property of a finite set of forbidden subgraphs 
that imply the existence of a dominating cycle.

\begin{lem}
\label{tree lem}
Let $\calH$ be a finite set of connected graphs, 
and suppose 
that there exists a positive integer $n = n(\calH)$ 
such that every $2$-connected $\calH$-free graph of order at least $n$ has a dominating cycle. 
\vspace{-5pt}
\begin{enumerate}[{\upshape(i)}]
\item Then $\calH$ contains a tree $T$ with $\Delta (T) \leq 3$ and $|V_{3}(T)| \leq 1$.
\vspace{-5pt}
\item If $|\calH|=2$ and $\calH$ contains a graph with diameter at least $3$, then the other one 
is an induced subgraph of $K_{1}+3K_{2}$, 
or 
isomorphic to $K_{4}^{-}$. 
\end{enumerate}
\end{lem}

\begin{figure}
\begin{center}
\unitlength 0.1in
\begin{picture}( 50.9500,  8.5000)( 16.0000,-19.7000)
%
\special{pn 8}%
\special{sh 1.000}%
\special{ar 2020 1636 40 40  0.0000000 6.2831853}%
%
\special{pn 8}%
\special{sh 1.000}%
\special{ar 2220 1636 40 40  0.0000000 6.2831853}%
%
\special{pn 8}%
\special{sh 1.000}%
\special{ar 2220 1386 40 40  0.0000000 6.2831853}%
%
\special{pn 8}%
\special{sh 1.000}%
\special{ar 2220 1886 40 40  0.0000000 6.2831853}%
%
\special{pn 8}%
\special{sh 1.000}%
\special{ar 2420 1386 40 40  0.0000000 6.2831853}%
%
\special{pn 8}%
\special{sh 1.000}%
\special{ar 2420 1886 40 40  0.0000000 6.2831853}%
%
\special{pn 8}%
\special{sh 1.000}%
\special{ar 2420 1636 40 40  0.0000000 6.2831853}%
%
\special{pn 8}%
\special{sh 1.000}%
\special{ar 2820 1636 40 40  0.0000000 6.2831853}%
%
\special{pn 8}%
\special{sh 1.000}%
\special{ar 2820 1386 40 40  0.0000000 6.2831853}%
%
\special{pn 8}%
\special{sh 1.000}%
\special{ar 2820 1886 40 40  0.0000000 6.2831853}%
%
\special{pn 8}%
\special{sh 1.000}%
\special{ar 3020 1636 40 40  0.0000000 6.2831853}%
%
\special{pn 8}%
\special{pa 3020 1636}%
\special{pa 2820 1636}%
\special{fp}%
\special{pa 2820 1636}%
\special{pa 2820 1636}%
\special{fp}%
%
\special{pn 8}%
\special{pa 2020 1636}%
\special{pa 2420 1636}%
\special{fp}%
%
\special{pn 8}%
\special{pa 2020 1636}%
\special{pa 2220 1386}%
\special{fp}%
%
\special{pn 8}%
\special{pa 2220 1886}%
\special{pa 2020 1636}%
\special{fp}%
%
\special{pn 8}%
\special{pa 2220 1386}%
\special{pa 2420 1386}%
\special{fp}%
%
\special{pn 8}%
\special{pa 2220 1886}%
\special{pa 2420 1886}%
\special{fp}%
%
\special{pn 8}%
\special{pa 2820 1886}%
\special{pa 3020 1636}%
\special{fp}%
\special{pa 3020 1636}%
\special{pa 2820 1386}%
\special{fp}%
%
\special{pn 8}%
\special{pa 2820 1386}%
\special{pa 2720 1386}%
\special{fp}%
%
\special{pn 8}%
\special{pa 2420 1386}%
\special{pa 2510 1386}%
\special{fp}%
%
\special{pn 8}%
\special{pa 2520 1386}%
\special{pa 2720 1386}%
\special{dt 0.045}%
%
\special{pn 8}%
\special{pa 2820 1636}%
\special{pa 2720 1636}%
\special{fp}%
%
\special{pn 8}%
\special{pa 2420 1636}%
\special{pa 2510 1636}%
\special{fp}%
%
\special{pn 8}%
\special{pa 2520 1636}%
\special{pa 2720 1636}%
\special{dt 0.045}%
%
\special{pn 8}%
\special{pa 2820 1886}%
\special{pa 2720 1886}%
\special{fp}%
%
\special{pn 8}%
\special{pa 2420 1886}%
\special{pa 2510 1886}%
\special{fp}%
%
\special{pn 8}%
\special{pa 2520 1886}%
\special{pa 2720 1886}%
\special{dt 0.045}%
\put(25.3000,-12.0500){\makebox(0,0){{\scriptsize $s$ vertices}}}%
%
\special{pn 8}%
\special{pa 2610 1290}%
\special{pa 2850 1290}%
\special{fp}%
\special{sh 1}%
\special{pa 2850 1290}%
\special{pa 2784 1270}%
\special{pa 2798 1290}%
\special{pa 2784 1310}%
\special{pa 2850 1290}%
\special{fp}%
%
\special{pn 8}%
\special{pa 2650 1290}%
\special{pa 2190 1290}%
\special{fp}%
\special{sh 1}%
\special{pa 2190 1290}%
\special{pa 2258 1310}%
\special{pa 2244 1290}%
\special{pa 2258 1270}%
\special{pa 2190 1290}%
\special{fp}%
\put(25.4000,-20.5500){\makebox(0,0){{\small $A_{s}$}}}%
%
\special{pn 8}%
\special{sh 1.000}%
\special{ar 3690 1390 40 40  0.0000000 6.2831853}%
%
\special{pn 8}%
\special{sh 1.000}%
\special{ar 3890 1390 40 40  0.0000000 6.2831853}%
%
\special{pn 8}%
\special{pa 3890 1390}%
\special{pa 3690 1390}%
\special{fp}%
%
\special{pn 8}%
\special{sh 1.000}%
\special{ar 4090 1390 40 40  0.0000000 6.2831853}%
%
\special{pn 8}%
\special{sh 1.000}%
\special{ar 4290 1390 40 40  0.0000000 6.2831853}%
%
\special{pn 8}%
\special{pa 4290 1390}%
\special{pa 4090 1390}%
\special{fp}%
%
\special{pn 8}%
\special{sh 1.000}%
\special{ar 4690 1390 40 40  0.0000000 6.2831853}%
%
\special{pn 8}%
\special{sh 1.000}%
\special{ar 4890 1390 40 40  0.0000000 6.2831853}%
%
\special{pn 8}%
\special{pa 4890 1390}%
\special{pa 4690 1390}%
\special{fp}%
%
\special{pn 8}%
\special{sh 1.000}%
\special{ar 4090 1890 40 40  0.0000000 6.2831853}%
%
\special{pn 8}%
\special{sh 1.000}%
\special{ar 4490 1890 40 40  0.0000000 6.2831853}%
%
\special{pn 8}%
\special{pa 4090 1890}%
\special{pa 3690 1390}%
\special{fp}%
\special{pa 3890 1390}%
\special{pa 4090 1890}%
\special{fp}%
\special{pa 4090 1890}%
\special{pa 4090 1390}%
\special{fp}%
\special{pa 4290 1390}%
\special{pa 4090 1890}%
\special{fp}%
\special{pa 4090 1890}%
\special{pa 4690 1390}%
\special{fp}%
\special{pa 4890 1390}%
\special{pa 4090 1890}%
\special{fp}%
%
\special{pn 8}%
\special{pa 3690 1390}%
\special{pa 4490 1890}%
\special{fp}%
\special{pa 4490 1890}%
\special{pa 3890 1390}%
\special{fp}%
\special{pa 4090 1390}%
\special{pa 4490 1900}%
\special{fp}%
\special{pa 4490 1900}%
\special{pa 4290 1390}%
\special{fp}%
\special{pa 4690 1390}%
\special{pa 4490 1890}%
\special{fp}%
\special{pa 4490 1890}%
\special{pa 4890 1390}%
\special{fp}%
%
\special{pn 10}%
\special{pa 4400 1390}%
\special{pa 4580 1390}%
\special{dt 0.045}%
%
\special{pn 8}%
\special{pa 4380 1300}%
\special{pa 4930 1300}%
\special{fp}%
\special{sh 1}%
\special{pa 4930 1300}%
\special{pa 4864 1280}%
\special{pa 4878 1300}%
\special{pa 4864 1320}%
\special{pa 4930 1300}%
\special{fp}%
%
\special{pn 8}%
\special{pa 4540 1300}%
\special{pa 3660 1300}%
\special{fp}%
\special{sh 1}%
\special{pa 3660 1300}%
\special{pa 3728 1320}%
\special{pa 3714 1300}%
\special{pa 3728 1280}%
\special{pa 3660 1300}%
\special{fp}%
\put(42.9000,-12.1000){\makebox(0,0){{\scriptsize $sK_{2}$}}}%
%
\special{pn 8}%
\special{sh 1.000}%
\special{ar 5616 1896 40 40  0.0000000 6.2831853}%
%
\special{pn 8}%
\special{sh 1.000}%
\special{ar 5816 1896 40 40  0.0000000 6.2831853}%
%
\special{pn 8}%
\special{pa 5816 1896}%
\special{pa 5616 1896}%
\special{fp}%
%
\special{pn 8}%
\special{sh 1.000}%
\special{ar 6216 1896 40 40  0.0000000 6.2831853}%
%
\special{pn 8}%
\special{sh 1.000}%
\special{ar 6416 1896 40 40  0.0000000 6.2831853}%
%
\special{pn 8}%
\special{pa 6416 1896}%
\special{pa 6216 1896}%
\special{fp}%
%
\special{pn 8}%
\special{sh 1.000}%
\special{ar 5810 1636 40 40  0.0000000 6.2831853}%
%
\special{pn 8}%
\special{sh 1.000}%
\special{ar 6210 1636 40 40  0.0000000 6.2831853}%
\put(42.9000,-12.1000){\makebox(0,0){{\scriptsize $sK_{2}$}}}%
%
\special{pn 8}%
\special{pa 5810 1636}%
\special{pa 5730 1346}%
\special{fp}%
%
\special{pn 8}%
\special{pa 5810 1636}%
\special{pa 6190 1356}%
\special{fp}%
%
\special{pn 8}%
\special{pa 6210 1636}%
\special{pa 5810 1376}%
\special{fp}%
%
\special{pn 8}%
\special{pa 6210 1636}%
\special{pa 6270 1356}%
\special{fp}%
%
\special{pn 8}%
\special{sh 0}%
\special{ar 6000 1286 340 140  0.0000000 6.2831853}%
%
\put(59.9000,-12.8500){\makebox(0,0){}}%
\put(42.9000,-12.1000){\makebox(0,0){{\scriptsize $sK_{2}$}}}%
\put(60.0500,-12.8500){\makebox(0,0){{\scriptsize $K_{s}$}}}%
%
\special{pn 8}%
\special{pa 5816 1636}%
\special{pa 5616 1896}%
\special{fp}%
\special{pa 5816 1896}%
\special{pa 5816 1636}%
\special{fp}%
\special{pa 5816 1636}%
\special{pa 6216 1896}%
\special{fp}%
\special{pa 6416 1896}%
\special{pa 5816 1636}%
\special{fp}%
%
\special{pn 8}%
\special{pa 6216 1636}%
\special{pa 6416 1896}%
\special{fp}%
\special{pa 6216 1896}%
\special{pa 6216 1636}%
\special{fp}%
\special{pa 6216 1636}%
\special{pa 5816 1896}%
\special{fp}%
\special{pa 5616 1896}%
\special{pa 6216 1636}%
\special{fp}%
\put(58.6500,-15.0500){\makebox(0,0){{\scriptsize $+$}}}%
\put(61.5500,-15.0500){\makebox(0,0){{\scriptsize $+$}}}%
\put(42.9000,-20.5500){\makebox(0,0){{\small $A_{s}'$}}}%
\put(60.3000,-20.5500){\makebox(0,0){{\small $A_{s}''$}}}%
%
\special{pn 8}%
\special{pa 5810 1640}%
\special{pa 6210 1640}%
\special{fp}%
\end{picture}%
\caption{The graphs $A_{s}$, $A_{s}'$ and $A_{s}''$}
\label{As-graphs}
\end{center}
\end{figure}
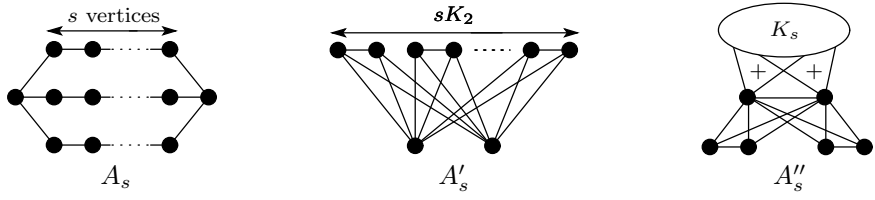

In order to prove Lemma~\ref{tree lem}, 
we define the following graphs $A_{s}$, $A_{s}'$ and $A_{s}''$ (see Figure~\ref{As-graphs}). 
Note that  each of $A_{s}$, $A_{s}'$ and $A_{s}''$ is $2$-connected and contains no dominating cycle. 
\begin{itemize}
\item
For each $s\geq 2$, 
let $A_{s}$ be the graph 
consisting of the union of three internally disjoint $P_{s+2}$'s that have the same two distinct end vertices. 
\vspace{-5pt}
\item
For each $s\geq 3$, let $A'_{s}=2K_{1}+sK_{2}$.
\vspace{-5pt}
\item
For each $s\geq 2$, let $A''_{s}=K_{2}+(2K_{2}\cup K_{s})$.
\end{itemize}

\medskip
\noindent
\textbf{Proof of Lemma~\ref{tree lem}.}~\rm{(i)} 
Let $m = \max\{|V(H)| : H \in \calH \}$, 
and let $n_{1} = \max \{n, m\}$. 
Since $A_{n_{1}}$ is a $2$-connected graph of order at least $n_{1} \ ( \ \ge n)$ having no dominating cycle, 
it follows that there exists a graph $H$ in $\calH$ 
such that $H \prec A_{n_{1}}$. 
Observe that, by the definition of $A_{n_{1}}$, all cycles in $A_{n_{1}}$ have $2n_{1} + 2 \ ( \ \ge 2m + 2)$ vertices 
and the distance of two vertices with degree $3$ in $A_{n_{1}}$ is $n_{1}+1 \ ( \ \ge m + 1)$. 
These facts imply that $H$ is a tree with $\Delta (H)\leq 3$ and $|V_{3}(H)|\leq 1$.

\rm{(ii)} Write $\calH = \{H_{1}, H_{2}\}$, 
and assume that $\diam(H_{1}) \geq 3$. 
Let $n_{2} = \max \{n, 3\}$. 
Since $\diam(H_{1}) \geq 3$, 
we have $P_{4} \prec H_{1}$, 
and hence $A'_{n_{2}}$ does not contain $H_{1}$ as an induced subgraph 
because $A'_{n_{2}}$ contains no $P_{4}$ as an induced subgraph. 
Similarly, we see that $A''_{n_{2}}$ does not contain $H_{1}$ as an induced subgraph.
On the other hand, both of 
$A'_{n_{2}}$ and 
$A''_{n_{2}}$ are $2$-connected graphs of order at least $n_{2} \ ( \ \ge n)$ having no dominating cycle.
This implies that $H_{2}$ is a common induced subgraph of $A'_{n_{2}}$ and $A''_{n_{2}}$. 
Hence it is easy to check that $H_{2} \prec K_{1}+3K_{2}$ or $H_{2} \cong K_{4}^{-}$.
\qed

\begin{figure}
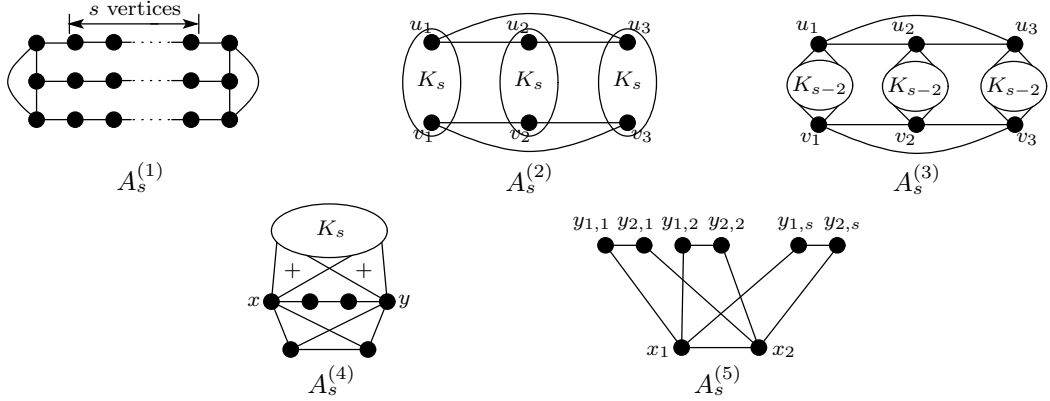

\begin{center}
\unitlength 0.1in
%
\caption{The graph $A_{s}^{(j)}$}
\label{As(j)-graphs}
\end{center}
\end{figure}

We further define five graphs of 2-connected graphs having no dominating cycle as follows (see Figure~\ref{As(j)-graphs}). 
\begin{itemize}
\item
For each $s \geq 2$, 
let $A^{(1)}_{s}$ be the graph which consists of two vertex-disjoint triangles 
connected by three vertex-disjoint paths of orders $s+2$, respectively. 
\vspace{-5pt}
\item
For each $s\geq 4$, let $G_{i}~(1\leq i\leq 3)$ be a complete graph of order $s$,
and let $A^{(2)}_{s}$ be the graph obtained from $G_{1}\cup G_{2}\cup G_{3}$
by joining $u_{i}$ to $u_{i+1}$ and $v_{i}$ to $v_{i+1}$ for $1\leq i\leq 3$,
where $u_{i}$ and $v_{i}$ are distinct two vertices of $G_{i}$ and $u_{4} = u_{1},v_{4} = v_{1}$.
\vspace{-5pt}
\item
For each $s \geq 4$, let $A^{(3)}_{s} = A^{(2)}_{s} - \{u_{i}v_{i} : 1\leq i\leq 3\}$. 
\vspace{-5pt}
\item
For each $s\geq 2$, 
let $A^{(4)}_{s}$ be the graph obtained from $K_{2} + (K_{2} \cup K_{s})$ 
by subdividing the edge $xy$ twice, where $\{x,y\}$ is the unique $2$-cut set of $K_{2} + (K_{2} \cup K_{s})$. 
\vspace{-5pt}
\item
For each $s \geq 3$, 
let $A^{(5)}_{s}$ be the graph defined by 
$V(A^{(5)}_{s})=\{x_{i},y_{i,j} : 1\leq i\leq 2,1\leq j\leq s\}$ 
and $E(A^{(5)}_{s})=\{x_{1}x_{2}\}\cup \{y_{1,j}y_{2,j} : 1\leq j\leq s\}\cup \{x_{i}y_{i,j} : 1\leq i\leq 2,1\leq j\leq s\}$.
\end{itemize}

By the definition of $A_{s}^{(j)}$, 
we can obtain the following lemma. 
(Since the proof is easy, we omit it.)

\begin{lem}
\label{A_{s}^{(1)}, A_{s}^{(2)} and A_{s}^{(3)}}
\begin{enumerate}[{\upshape(i)}]
\item $A_{s}^{(1)}$, $A_{s}^{(2)}$ and $A_{s}^{(3)}$ are $K_{1, 3}$-free graphs. 
Furthermore, $A_{s}^{(1)}$ is $K_{4}^{-}$-free. 
\vspace{-5pt}
\item Every connected induced subgraph of $A_{s}^{(1)}$ with at most $s$ vertices is also an induced subgraph of $N_{i,j,k}$ for some integers $i$, $j$ and $k$.
In particular, every induced subtree of $A_{s}^{(1)}$ with at most $s$ vertices is a path.
\vspace{-5pt}
\item Every induced subpath of $A_{s}^{(2)}$ has at most $6$ vertices. 
\vspace{-5pt}
\item $A_{s}^{(2)}$ contains neither $N_{1, 1, 2}$ nor $B_{1,3}$ as an induced subgraph.
\vspace{-5pt}
\item $A_{s}^{(3)}$ contains no $B_{2, 2}$ as an induced subgraph. 
\vspace{-5pt}
\item $A_{s}^{(4)}$ is $P_{5}$-free and $A_{s}^{(5)}$ is $K_{3}$-free. 
\end{enumerate}
\end{lem}

Now we are ready to prove Theorem~\ref{2-conn (necessity)}. 

\medskip
\noindent
\textbf{Proof of Theorem~\ref{2-conn (necessity)}.}~Let $\cal{H}$ be a set of two connected graphs, 
and suppose that 
there exists a positive integer $n_{0} = n_{0}(\calH)$ such that 
every $2$-connected $\cal{H}$-free graph of order at least $n_{0}$ has a dominating cycle. 
We show that 
$\calH \le \calH_{i}$ for some $i$ with $1 \le i \le 7$.

If $\calH$ contains $P_{3}$, 
then $\calH \leq \calH_{1}$. 
Thus we may assume that $\calH$ does not contain $P_{3}$. 
Write $\calH = \{H_{1}, H_{2}\}$, 
and 
let $n = \max \{n_{0}, 4, |V(H_{1})|, |V(H_{2})|\}$. 
Then 
for each $j$ with $1 \leq j \leq 5$, 
\begin{align}
\label{s_{j}}
\textup{$H_{1} \prec A_{n}^{(j)}$ or $H_{2} \prec A_{n}^{(j)}$}
\end{align}
because $A_{n}^{(j)}$ is a 2-connected graph of order at least $n \ ( \ \ge n_{0})$ having no dominating cycle.

We divide the proof into two cases according as 
$\calH$ contains $K_{1, 3}$ or not.

\medskip
\noindent
\textbf{Case 1.} $H_{i}$ is isomorphic to $K_{1,3}$ for some $i$ with $i \in \{1, 2\}$.

We may assume that $H_{1} \cong K_{1, 3}$. 
Then it follows from Lemma~\ref{A_{s}^{(1)}, A_{s}^{(2)} and A_{s}^{(3)}} (i) 
and 
(\ref{s_{j}}) 
that 
$H_{2} \prec A_{n}^{(j)}$ for each $1 \le j \le 3$, 
i.e., 
$H_{2}$ is a common induced subgraph 
of $A_{n}^{(1)}$, $A_{n}^{(2)}$ and $A_{n}^{(3)}$.
Since $n \ge |V(H_{2})|$, 
it follows from Lemma~\ref{A_{s}^{(1)}, A_{s}^{(2)} and A_{s}^{(3)}} (ii)
that $H_{2}$ is an induced subgraph of $N_{i, j, k}$ for some integers $i$, $j$ and $k$. 
This together with 
Lemma~\ref{A_{s}^{(1)}, A_{s}^{(2)} and A_{s}^{(3)}} (iii) implies that 
$H_{2}$ is an induced subgraph of either $Z_{4}$, $B_{1,3}$ or $N_{2,2,2}$.
If $H_{2} \prec Z_{4}$, then $\calH \leq \calH_{1}$; 
if $H_{2}\prec B_{1,3}$, then by Lemma~\ref{A_{s}^{(1)}, A_{s}^{(2)} and A_{s}^{(3)}} (iv),
either $H_{2} \prec P_{6}$, $H_{2} \prec Z_{3}$ or $H_{2} \prec B_{1,2}$,
and hence $\calH\leq \calH_{1}$ or $\calH\leq \calH_{2}$;
if $H_{2 }\prec N_{2,2,2}$, then by Lemma~\ref{A_{s}^{(1)}, A_{s}^{(2)} and A_{s}^{(3)}} (iv) and (v),
either $H_{2} \prec P_{6}$, $H_{2} \prec Z_{2}$, $H_{2} \prec B_{1,2}$ or $H_{2} \prec N_{1,1,1}$,
and hence $\calH\leq \calH_{1}$, $\calH\leq \calH_{2}$ or $\calH\leq \calH_{3}$.

\medskip
\noindent
\textbf{Case 2.} $H_{i}$ is not isomorphic to $K_{1,3}$ for $i \in \{1, 2\}$.

By Lemma~\ref{tree lem} (i), 
we may assume that $H_{1}$ is a tree with $\Delta (H_{1})\leq 3$ and $|V_{3}(H_{1})|\leq 1$.
Since $H_{1} \not\cong P_{3}$ and $H_{1} \not\cong K_{1, 3}$ by the assumption of Case 2,
we have $\mbox{diam}(H_{1})\geq 3$.
Hence by Lemma~\ref{tree lem} (ii), 
$H_{2}$ 
is 
an induced subgraph of $K_{1}+3K_{2}$ 
or 
isomorphic to $K_{4}^{-}$. 
In particular, 
this implies that $H_{2}$ has a triangle 
(note that if $H_{2}$ is a tree, then either $H_{2} \cong K_{1,3}$ or $H_{2} \cong P_{3}$, a contradiction).
Therefore, 
it 
follows from Lemma~\ref{A_{s}^{(1)}, A_{s}^{(2)} and A_{s}^{(3)}} (vi) and (\ref{s_{j}}) that $H_{1} \prec A_{n}^{(5)}$. 
Combining this with the fact that 
$\Delta (H_{1})\leq 3$ and $|V_{3}(H_{1})|\leq 1$, 
we have $H_{1} \prec K^{*}_{1,3}$.

We divide the proof of Case 2 into three cases according as 
$H_{2} \prec K_{1}+3K_{2}$ and $H_{1}$ is not a path, 
$H_{2} \prec K_{1}+3K_{2}$ and $H_{1}$ is a path, 
or 
$H_{2} \cong K_{4}^{-}$.

\medskip
\noindent
\textbf{Subcase 2.1.} $H_{2}$ is an induced subgraph of $K_{1}+3K_{2}$ and $H_{1}$ is not a path.

Since $H_{1}$ is a tree which is not a path,
$H_{1}$ is not an induced subgraph of $A_{n}^{(1)}$ by Lemma~\ref{A_{s}^{(1)}, A_{s}^{(2)} and A_{s}^{(3)}} (ii).
This together with (\ref{s_{j}}) implies that $H_{2} \prec A^{(1)}_{n}$.
Combining this with the assumption that $H_{2}$ is an induced subgraph of $K_{1}+3K_{2}$, 
we see that $H_{2} \prec Z_{1}$. 
Since $H_{1} \prec K_{1, 3}^{*}$, 
$\calH \leq \calH_{5}$.

\medskip
\noindent
\textbf{Subcase 2.2.} $H_{2}$ is an induced subgraph of $K_{1}+3K_{2}$ and $H_{1}$ is a path.

If $H_{1}$ is a path of order at most $4$, then $\calH \leq \calH_{4}$ (note that $K_{1}+3K_{2} \cong W$).
Thus we may assume that $H_{1}$ is a path of order $5$ because $H_{1} \prec K^{*}_{1,3}$. 
Then by Lemma~\ref{A_{s}^{(1)}, A_{s}^{(2)} and A_{s}^{(3)}} (vi) and 
(\ref{s_{j}}), $H_{2} \prec A^{(4)}_{n}$.
Since $H_{2}\prec K_{1}+3K_{2}$ and $A_{n}^{(4)}$ contains no $K_{1}+3K_{2}$ as an induced subgraph,
$H_{2}$ is an induced subgraph of $K_{1}+(K_{1}\cup 2K_{2})$.
Thus $\calH \leq \calH_{6}$ (note that $K_{1}+(K_{1}\cup 2K_{2}) \cong W^{*}$).

\medskip
\noindent
\textbf{Subcase 2.3.} $H_{2}$ is isomorphic to $K_{4}^{-}$.

Then by Lemma~\ref{A_{s}^{(1)}, A_{s}^{(2)} and A_{s}^{(3)}} (i), 
(\ref{s_{j}}) and the assumption of Subcase 2.3, 
we have $H_{1} \prec A^{(1)}_{n}$.
Since $H_{1} \prec K^{*}_{1,3}$, $H_{1}$ is a path of order at most $5$ by Lemma~\ref{A_{s}^{(1)}, A_{s}^{(2)} and A_{s}^{(3)}} (ii).
Thus $\calH \leq \calH_{7}$.

This completes the proof of Theorem~\ref{2-conn (necessity)}.
\qed

\section{Proof of Theorem~\ref{2-conn (sufficiency)}}
\label{proof of sufficiency}

In this section, 
we prove Theorem~\ref{2-conn (sufficiency)}. 
To prove this, 
we show that the following theorems hold, 
which immediately imply Theorem~\ref{2-conn (sufficiency)} 
(note that we actually prove slightly stronger statements.) 

\begin{thm}
\label{K_{1,3}-free}
If $\calH \in \big\{ \calH_{i} : 1 \le i \le 3 \big\}$, 
then 
every longest cycle of a $2$-connected $\calH$-free graph is a dominating cycle of the graph. 
\end{thm}

\begin{thm}
\label{P_{4}, W, long}
Every longest cycle of a $2$-connected $\{ P_{4}, W \}$-free graph is a dominating cycle of the graph. 
\end{thm}

\begin{thm}
\label{K_{1,3}^{**}, Z_{1}}
A longest cycle of a $2$-connected $\{ K_{1, 3}^{**}, Z_{1} \}$-free graph is a dominating cycle of the graph.
\end{thm}

We will prove Theorems~\ref{K_{1,3}-free} and \ref{P_{4}, W, long} 
in Subsections~\ref{K_{1, 3}-free graphs} and \ref{Proof of P_{4}, W}, respectively, 
and we will prove Theorem~\ref{K_{1,3}^{**}, Z_{1}} in Subsections~\ref{Proof of K_{1, 3}^{**}, Z_{1}}--\ref{Proof of tri th}.

\subsection{$K_{1, 3}$-free graphs}
\label{K_{1, 3}-free graphs}

In this subsection, 
we prove Theorem~\ref{K_{1,3}-free}. 
In order to prove it, 
we use some concepts and known results.

In \cite{R}, Ryj\'{a}\v{c}ek introduced the concept of a closure for claw-free graphs as follows. 
Let $G$ be a claw-free graph. 
For each vertex $v$ of $G$, 
$G[N_{G}(v)]$ has at most two components; 
otherwise $G$ contains a $K_{1, 3}$ as an induced subgraph. 
If $G[N_{G}(v)]$ has two components, both of them must be cliques. 
In the case that $G[N_{G}(v)]$ is connected, 
we add edges joining all pairs of nonadjacent vertices in $N_{G}(v)$. 
The \textit{closure} $\cl(G)$ of $G$ is a graph obtained
by recursively repeating this operation, as long as this is possible.
In \cite{R}, 
it was shown that the closure of a graph has the following property.

\begin{Thm}[Ryj\'{a}\v{c}ek~\cite{R}] 
\label{R-closure}
If $G$ is a claw-free graph, then the following hold. 
\vspace{-5pt}
\begin{enumerate}[{\upshape(i)}]
\item $\cl(G)$ is well-defined, $($i.e., uniquely defined$)$.
\vspace{-5pt} 
\item $c(G) = c(\cl(G))$. 
\end{enumerate}
\end{Thm}

\begin{figure}
\begin{center}
\unitlength 0.1in
\begin{picture}( 46.7500,  9.7500)( 16.0000,-20.6500)
%
\special{pn 8}%
\special{sh 1.000}%
\special{ar 4220 1840 40 40  0.0000000 6.2831853}%
%
\special{pn 8}%
\special{sh 1.000}%
\special{ar 4220 1690 40 40  0.0000000 6.2831853}%
%
\special{pn 8}%
\special{sh 1.000}%
\special{ar 4220 1540 40 40  0.0000000 6.2831853}%
%
\special{pn 8}%
\special{sh 1.000}%
\special{ar 4220 1390 40 40  0.0000000 6.2831853}%
%
\special{pn 8}%
\special{sh 1.000}%
\special{ar 4720 1840 40 40  0.0000000 6.2831853}%
%
\special{pn 8}%
\special{sh 1.000}%
\special{ar 4720 1390 40 40  0.0000000 6.2831853}%
%
\special{pn 8}%
\special{sh 1.000}%
\special{ar 4470 1620 40 40  0.0000000 6.2831853}%
%
\special{pn 8}%
\special{sh 1.000}%
\special{ar 4720 1620 40 40  0.0000000 6.2831853}%
%
\special{pn 8}%
\special{sh 1.000}%
\special{ar 4470 1240 40 40  0.0000000 6.2831853}%
%
\special{pn 8}%
\special{sh 1.000}%
\special{ar 4470 1990 40 40  0.0000000 6.2831853}%
%
\special{pn 8}%
\special{sh 1.000}%
\special{ar 3920 1840 40 40  0.0000000 6.2831853}%
%
\special{pn 8}%
\special{sh 1.000}%
\special{ar 3920 1390 40 40  0.0000000 6.2831853}%
%
\special{pn 8}%
\special{sh 1.000}%
\special{ar 3670 1620 40 40  0.0000000 6.2831853}%
%
\special{pn 8}%
\special{sh 1.000}%
\special{ar 3920 1620 40 40  0.0000000 6.2831853}%
%
\special{pn 8}%
\special{sh 1.000}%
\special{ar 3670 1240 40 40  0.0000000 6.2831853}%
%
\special{pn 8}%
\special{sh 1.000}%
\special{ar 3670 1990 40 40  0.0000000 6.2831853}%
%
\special{pn 8}%
\special{sh 1.000}%
\special{ar 3420 1840 40 40  0.0000000 6.2831853}%
%
\special{pn 8}%
\special{sh 1.000}%
\special{ar 3420 1390 40 40  0.0000000 6.2831853}%
%
\special{pn 8}%
\special{sh 1.000}%
\special{ar 3420 1620 40 40  0.0000000 6.2831853}%
%
\special{pn 8}%
\special{pa 3420 1390}%
\special{pa 3920 1390}%
\special{fp}%
\special{pa 3920 1390}%
\special{pa 3920 1840}%
\special{fp}%
\special{pa 3920 1840}%
\special{pa 3420 1840}%
\special{fp}%
\special{pa 3420 1840}%
\special{pa 3420 1390}%
\special{fp}%
\special{pa 3420 1390}%
\special{pa 3670 1240}%
\special{fp}%
\special{pa 3670 1240}%
\special{pa 3920 1390}%
\special{fp}%
\special{pa 3920 1840}%
\special{pa 3670 1990}%
\special{fp}%
\special{pa 3670 1990}%
\special{pa 3420 1840}%
\special{fp}%
\special{pa 3670 1240}%
\special{pa 3670 1990}%
\special{fp}%
%
\special{pn 8}%
\special{pa 4220 1390}%
\special{pa 4720 1390}%
\special{fp}%
\special{pa 4720 1390}%
\special{pa 4720 1840}%
\special{fp}%
\special{pa 4720 1840}%
\special{pa 4220 1840}%
\special{fp}%
\special{pa 4220 1840}%
\special{pa 4220 1390}%
\special{fp}%
\special{pa 4220 1390}%
\special{pa 4470 1240}%
\special{fp}%
\special{pa 4470 1240}%
\special{pa 4720 1390}%
\special{fp}%
\special{pa 4720 1840}%
\special{pa 4470 1990}%
\special{fp}%
\special{pa 4470 1990}%
\special{pa 4220 1840}%
\special{fp}%
\special{pa 4470 1240}%
\special{pa 4470 1990}%
\special{fp}%
%
\special{pn 8}%
\special{sh 1.000}%
\special{ar 3120 1840 40 40  0.0000000 6.2831853}%
%
\special{pn 8}%
\special{sh 1.000}%
\special{ar 3120 1390 40 40  0.0000000 6.2831853}%
%
\special{pn 8}%
\special{sh 1.000}%
\special{ar 2870 1620 40 40  0.0000000 6.2831853}%
%
\special{pn 8}%
\special{sh 1.000}%
\special{ar 3120 1620 40 40  0.0000000 6.2831853}%
%
\special{pn 8}%
\special{sh 1.000}%
\special{ar 2870 1240 40 40  0.0000000 6.2831853}%
%
\special{pn 8}%
\special{sh 1.000}%
\special{ar 2870 1990 40 40  0.0000000 6.2831853}%
%
\special{pn 8}%
\special{sh 1.000}%
\special{ar 2620 1840 40 40  0.0000000 6.2831853}%
%
\special{pn 8}%
\special{sh 1.000}%
\special{ar 2620 1390 40 40  0.0000000 6.2831853}%
%
\special{pn 8}%
\special{sh 1.000}%
\special{ar 2620 1620 40 40  0.0000000 6.2831853}%
%
\special{pn 8}%
\special{pa 2620 1390}%
\special{pa 3120 1390}%
\special{fp}%
\special{pa 3120 1390}%
\special{pa 3120 1840}%
\special{fp}%
\special{pa 3120 1840}%
\special{pa 2620 1840}%
\special{fp}%
\special{pa 2620 1840}%
\special{pa 2620 1390}%
\special{fp}%
\special{pa 2620 1390}%
\special{pa 2870 1240}%
\special{fp}%
\special{pa 2870 1240}%
\special{pa 3120 1390}%
\special{fp}%
\special{pa 3120 1840}%
\special{pa 2870 1990}%
\special{fp}%
\special{pa 2870 1990}%
\special{pa 2620 1840}%
\special{fp}%
\special{pa 2870 1240}%
\special{pa 2870 1990}%
\special{fp}%
%
\special{pn 8}%
\special{sh 1.000}%
\special{ar 2320 1840 40 40  0.0000000 6.2831853}%
%
\special{pn 8}%
\special{sh 1.000}%
\special{ar 2320 1390 40 40  0.0000000 6.2831853}%
%
\special{pn 8}%
\special{sh 1.000}%
\special{ar 2070 1620 40 40  0.0000000 6.2831853}%
%
\special{pn 8}%
\special{sh 1.000}%
\special{ar 2320 1620 40 40  0.0000000 6.2831853}%
%
\special{pn 8}%
\special{sh 1.000}%
\special{ar 2070 1240 40 40  0.0000000 6.2831853}%
%
\special{pn 8}%
\special{sh 1.000}%
\special{ar 2070 1990 40 40  0.0000000 6.2831853}%
%
\special{pn 8}%
\special{sh 1.000}%
\special{ar 1820 1840 40 40  0.0000000 6.2831853}%
%
\special{pn 8}%
\special{sh 1.000}%
\special{ar 1820 1390 40 40  0.0000000 6.2831853}%
%
\special{pn 8}%
\special{sh 1.000}%
\special{ar 1820 1620 40 40  0.0000000 6.2831853}%
%
\special{pn 8}%
\special{pa 1820 1390}%
\special{pa 2320 1390}%
\special{fp}%
\special{pa 2320 1390}%
\special{pa 2320 1840}%
\special{fp}%
\special{pa 2320 1840}%
\special{pa 1820 1840}%
\special{fp}%
\special{pa 1820 1840}%
\special{pa 1820 1390}%
\special{fp}%
\special{pa 1820 1390}%
\special{pa 2070 1240}%
\special{fp}%
\special{pa 2070 1240}%
\special{pa 2320 1390}%
\special{fp}%
\special{pa 2320 1840}%
\special{pa 2070 1990}%
\special{fp}%
\special{pa 2070 1990}%
\special{pa 1820 1840}%
\special{fp}%
\special{pa 2070 1240}%
\special{pa 2070 1990}%
\special{fp}%
%
\special{pn 8}%
\special{pa 2070 1240}%
\special{pa 2058 1272}%
\special{pa 2046 1302}%
\special{pa 2034 1334}%
\special{pa 2022 1364}%
\special{pa 2012 1396}%
\special{pa 2002 1426}%
\special{pa 1994 1458}%
\special{pa 1986 1488}%
\special{pa 1980 1520}%
\special{pa 1976 1550}%
\special{pa 1972 1580}%
\special{pa 1970 1612}%
\special{pa 1972 1642}%
\special{pa 1974 1674}%
\special{pa 1978 1704}%
\special{pa 1984 1736}%
\special{pa 1992 1766}%
\special{pa 2000 1798}%
\special{pa 2010 1828}%
\special{pa 2020 1860}%
\special{pa 2032 1890}%
\special{pa 2042 1920}%
\special{pa 2054 1952}%
\special{pa 2068 1982}%
\special{pa 2070 1990}%
\special{sp}%
%
\special{pn 8}%
\special{pa 4470 1240}%
\special{pa 4458 1272}%
\special{pa 4446 1302}%
\special{pa 4434 1334}%
\special{pa 4422 1364}%
\special{pa 4412 1396}%
\special{pa 4402 1426}%
\special{pa 4394 1458}%
\special{pa 4386 1488}%
\special{pa 4380 1520}%
\special{pa 4376 1550}%
\special{pa 4372 1580}%
\special{pa 4370 1612}%
\special{pa 4372 1642}%
\special{pa 4374 1674}%
\special{pa 4378 1704}%
\special{pa 4384 1736}%
\special{pa 4392 1766}%
\special{pa 4400 1798}%
\special{pa 4410 1828}%
\special{pa 4420 1860}%
\special{pa 4432 1890}%
\special{pa 4442 1920}%
\special{pa 4454 1952}%
\special{pa 4468 1982}%
\special{pa 4470 1990}%
\special{sp}%
%
\special{pn 8}%
\special{pa 2320 1390}%
\special{pa 2302 1420}%
\special{pa 2282 1450}%
\special{pa 2266 1478}%
\special{pa 2250 1508}%
\special{pa 2238 1538}%
\special{pa 2228 1566}%
\special{pa 2222 1596}%
\special{pa 2220 1626}%
\special{pa 2224 1654}%
\special{pa 2232 1684}%
\special{pa 2244 1712}%
\special{pa 2258 1742}%
\special{pa 2276 1772}%
\special{pa 2294 1800}%
\special{pa 2312 1830}%
\special{pa 2320 1840}%
\special{sp}%
%
\special{pn 8}%
\special{pa 3120 1390}%
\special{pa 3102 1420}%
\special{pa 3082 1450}%
\special{pa 3066 1478}%
\special{pa 3050 1508}%
\special{pa 3038 1538}%
\special{pa 3028 1566}%
\special{pa 3022 1596}%
\special{pa 3020 1626}%
\special{pa 3024 1654}%
\special{pa 3032 1684}%
\special{pa 3044 1712}%
\special{pa 3058 1742}%
\special{pa 3076 1772}%
\special{pa 3094 1800}%
\special{pa 3112 1830}%
\special{pa 3120 1840}%
\special{sp}%
%
\special{pn 8}%
\special{pa 4720 1390}%
\special{pa 4702 1420}%
\special{pa 4682 1450}%
\special{pa 4666 1478}%
\special{pa 4650 1508}%
\special{pa 4638 1538}%
\special{pa 4628 1566}%
\special{pa 4622 1596}%
\special{pa 4620 1626}%
\special{pa 4624 1654}%
\special{pa 4632 1684}%
\special{pa 4644 1712}%
\special{pa 4658 1742}%
\special{pa 4676 1772}%
\special{pa 4694 1800}%
\special{pa 4712 1830}%
\special{pa 4720 1840}%
\special{sp}%
\put(20.7000,-21.5000){\makebox(0,0){{\small $F_{1}$}}}%
\put(28.7000,-21.5000){\makebox(0,0){{\small $F_{2}$}}}%
\put(36.7000,-21.5000){\makebox(0,0){{\small $F_{3}$}}}%
\put(44.7000,-21.5000){\makebox(0,0){{\small $F_{4}$}}}%
%
\special{pn 8}%
\special{ar 5536 1856 490 140  0.0000000 6.2831853}%
%
\special{pn 8}%
\special{ar 5536 1236 490 140  0.0000000 6.2831853}%
%
\special{pn 8}%
\special{sh 1.000}%
\special{ar 5246 1796 40 40  0.0000000 6.2831853}%
%
\special{pn 8}%
\special{sh 1.000}%
\special{ar 5446 1796 40 40  0.0000000 6.2831853}%
%
\special{pn 8}%
\special{sh 1.000}%
\special{ar 5846 1796 40 40  0.0000000 6.2831853}%
%
\special{pn 8}%
\special{sh 1.000}%
\special{ar 5246 1546 40 40  0.0000000 6.2831853}%
%
\special{pn 8}%
\special{sh 1.000}%
\special{ar 5446 1546 40 40  0.0000000 6.2831853}%
%
\special{pn 8}%
\special{sh 1.000}%
\special{ar 5846 1546 40 40  0.0000000 6.2831853}%
%
\special{pn 8}%
\special{sh 1.000}%
\special{ar 5246 1296 40 40  0.0000000 6.2831853}%
%
\special{pn 8}%
\special{sh 1.000}%
\special{ar 5446 1296 40 40  0.0000000 6.2831853}%
%
\special{pn 8}%
\special{sh 1.000}%
\special{ar 5846 1296 40 40  0.0000000 6.2831853}%
%
\special{pn 8}%
\special{pa 5246 1796}%
\special{pa 5246 1296}%
\special{fp}%
\special{pa 5446 1296}%
\special{pa 5446 1796}%
\special{fp}%
\special{pa 5846 1796}%
\special{pa 5846 1296}%
\special{fp}%
%
\special{pn 8}%
\special{pa 5446 1296}%
\special{pa 5430 1326}%
\special{pa 5414 1356}%
\special{pa 5398 1386}%
\special{pa 5386 1416}%
\special{pa 5374 1446}%
\special{pa 5364 1476}%
\special{pa 5358 1506}%
\special{pa 5356 1536}%
\special{pa 5356 1566}%
\special{pa 5360 1596}%
\special{pa 5368 1626}%
\special{pa 5378 1656}%
\special{pa 5390 1686}%
\special{pa 5404 1718}%
\special{pa 5420 1748}%
\special{pa 5436 1778}%
\special{pa 5446 1796}%
\special{sp}%
%
\special{pn 8}%
\special{pa 5246 1296}%
\special{pa 5230 1326}%
\special{pa 5214 1356}%
\special{pa 5198 1386}%
\special{pa 5186 1416}%
\special{pa 5174 1446}%
\special{pa 5164 1476}%
\special{pa 5158 1506}%
\special{pa 5156 1536}%
\special{pa 5156 1566}%
\special{pa 5160 1596}%
\special{pa 5168 1626}%
\special{pa 5178 1656}%
\special{pa 5190 1686}%
\special{pa 5204 1718}%
\special{pa 5220 1748}%
\special{pa 5236 1778}%
\special{pa 5246 1796}%
\special{sp}%
%
\special{pn 8}%
\special{pa 5846 1296}%
\special{pa 5830 1326}%
\special{pa 5814 1356}%
\special{pa 5798 1386}%
\special{pa 5786 1416}%
\special{pa 5774 1446}%
\special{pa 5764 1476}%
\special{pa 5758 1506}%
\special{pa 5756 1536}%
\special{pa 5756 1566}%
\special{pa 5760 1596}%
\special{pa 5768 1626}%
\special{pa 5778 1656}%
\special{pa 5790 1686}%
\special{pa 5804 1718}%
\special{pa 5820 1748}%
\special{pa 5836 1778}%
\special{pa 5846 1796}%
\special{sp}%
%
\special{pn 10}%
\special{pa 5526 1546}%
\special{pa 5706 1546}%
\special{dt 0.045}%
\put(55.4500,-11.7500){\makebox(0,0){{\scriptsize $K_{s}$}}}%
\put(55.5000,-19.3000){\makebox(0,0){{\scriptsize $K_{s'}$}}}%
\put(62.7000,-15.1000){\makebox(0,0){{\scriptsize $2t+1$ cycles }}}%
\put(64.2000,-16.2000){\makebox(0,0){{\scriptsize of orders $3$}}}%
\put(55.5000,-21.5000){\makebox(0,0){{\small $F_{s,s',t}$}}}%
\end{picture}%
\caption{The graph $F_{i}$ and $F_{s, s', t}$}
\label{F-graphs}
\end{center}
\end{figure}
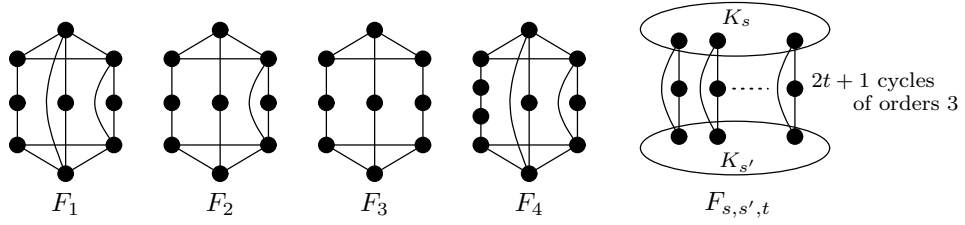

On the other hand, 
Brousek, Ryj\'{a}\v{c}ek and Favaron~\cite{BRF} 
characterized $2$-connected $\{K_{1,3},Z_{4}\}$-free graphs having no Hamilton cycle. 
Let $F_{1}$, $F_{2}$, $F_{3}$ and $F_{4}$ be the ones that are depicted in Figure~\ref{F-graphs}, 
and set $\F = \{F_{1}, F_{2}, F_{3}, F_{4}\}$. 
For each $s, s'$ and $t$ with $s' \ge s \ge 3$ and $1 \leq t \leq (s-1)/2$, 
let $F_{s, s', t}$ 
be the graph which consists of vertex-disjoint $K_{s}$ and $K_{s'}$ 
connected by $2t + 1$ vertex-disjoint cycles of orders $3$, respectively 
(see Figure~\ref{F-graphs}), 
and set 
$\F' = \{ F_{s,s',t} : s' \ge s \geq 3, 1 \leq t \leq (s-1)/2\}$.

\begin{Thm}[Brousek, Ryj\'{a}\v{c}ek and Favaron~\cite{BRF}]
\label{BRF}
Let $G$ be a $2$-connected $\{K_{1,3},Z_{4}\}$-free graph.
If $G$ is not Hamiltonian, 
then $G\in \F$ or $\cl(G) \in \F'$.
\end{Thm}

Now we prove Theorem~\ref{K_{1,3}-free}. 

\medskip
\noindent
\textbf{Proof of Theorem~\ref{K_{1,3}-free}.}~Let $\calH \in \{ \calH_{i} : 1 \le i \le 3\}$, 
and let $G$ be a $2$-connected $\calH$-free graph. 
We show that 
every longest cycle of $G$ is a dominating cycle of $G$. 
If $G$ is Hamiltonian, then the assertion clearly holds; 
thus 
we may assume that $G$ is not Hamiltonian. 
Then by Theorem~\ref{Be}, 
$\calH \neq \calH_{2} \ ( \ = \{K_{1, 3}, B_{1, 2}\})$ 
and 
$\calH \neq \calH_{3} \ ( \ = \{K_{1, 3}, N_{1, 1, 1}\})$. 
Thus $\calH = \calH_{1} \ ( \ = \{K_{1, 3}, Z_{4}\})$. 
Then it follows from Theorem~\ref{BRF} that $G \in \F$ or $\cl(G) \in \F'$. 
Since 
each graph $F$ in $\F \cup \F'$ has a longest cycle of order $|V(F)| - 1$, 
this together with Theorem~\ref{R-closure} (ii) implies that 
the longest cycle of $G$ has $|V(G)| - 1$ vertices; 
thus 
every longest cycle of $G$ is a dominating cycle. 
\qed

\subsection{$P_{4}$-free graphs}
\label{Proof of P_{4}, W}

In this subsection, we prove Theorem~\ref{P_{4}, W, long}.
To prove this, 
we use the following lemma concerning the property of $P_{4}$-free graphs.  
(In \cite{S}, a theorem which implies Lemma~\ref{P4-free} was proved by Seinseche, 
and also see \cite{EFFO, FFR}.)


\begin{Lem}
\label{P4-free}
Let $G$ be a $P_4$-free graph. 
If $G$ is $k$-connected and $|V(G)|\geq 2k$, 
then there exists a partition $\{A,B\}$ of $V(G)$ with $|A| \ge k$ and $|B| \geq k$ 
such that every vertex in $A$ is adjacent to each vertex in $B$.
\end{Lem}

Now we prove Theorem~\ref{P_{4}, W, long}.

\medskip
\noindent
\textbf{Proof of Theorem~\ref{P_{4}, W, long}.}~Let $G$ be a $2$-connected $\{P_{4}, W\}$-free graph. 
We may assume that $|V(G)| \ge 4$ (otherwise, the assertion clearly holds). 
Then by applying Lemma~\ref{P4-free} as $k=2$, 
there exists a partition $\{A,B\}$ of $V(G)$ with $|A| \ge 2$ and $|B| \geq 2$ 
such that every vertex in $A$ is adjacent to all vertices in $B$. 
By symmetry, we may assume that $|A|\geq |B|$.
Suppose that $G$ has a longest cycle $\ora{C}$ which is not a dominating cycle of $G$.

\begin{claim}
\label{cl:P4-1}
$B\subseteq V(C)$.
\end{claim}
\proof
Suppose that $B\not\subseteq V(C)$, and let $u\in B\setminus V(C)$.
Since $G[A'\cup B]$ is Hamiltonian for $A'\subseteq A$ with $|A'|=|B|$, we have $|V(C)|=c(G)\geq 2|B|$.
Hence $|V(C)\cap A| = |V(C) \setminus B| = |V(C) \setminus (B \setminus \{u\})| \geq |V(C)| - (|B|-1) 
\ge 2|B| - (|B| - 1) = |B| + 1$. 
Since $C-B \ ( \ = C - (B \setminus \{u\}))$ has at most $|B|-1$ components, 
this implies that 
there exists a vertex $x$ of $C$ such that $x, x^{+} \in A$.
Then the cycle $x^{+} \ora{C}xux^{+}$ is a longer cycle than $C$, a contradiction.
\qed

Since $C$ is not a dominating cycle of $G$, it follows from Claim~\ref{cl:P4-1} that 
there exist vertices $x_{1}, x_{2} \in A \setminus V(C)$ with $x_{1}x_{2} \in E(G)$. 
Moreover, by again Claim~\ref{cl:P4-1} and since $|B| \ge 2$, 
we can take distinct two vertices $u_{1}$ and $u_{2}$ in $B \cap V(C)$.

\begin{claim}
\label{cl:P4-2}
$N_{G}(\{x_{1}, x_{2}\}; \{u_{1}^{+}, u_{1}^{+2}, u_{2}^{+}, u_{2}^{+2}\}) = \emptyset$. 
In particular, 
$\{u_{1}^{+}, u_{1}^{+2}, u_{2}^{+}, u_{2}^{+2}\} \subseteq A$. 
\end{claim}
\proof
If there exists a vertex $u$ in $N_{G}(x_{1}; \{u_{1}^{+}, u_{1}^{+2}\})$, 
then $x_{1}u\ora{C}u_{1}x_{2}x_{1}$ is a cycle which contains $(V(C) \setminus \{u_{1}^{+}, u_{1}^{+2}\}) \cup \{u, x_{1}, x_{2}\}$. 
This contradicts the maximality of $|V(C)|$. 
Thus $N_{G}(x_{1}, \{u_{1}^{+}, u_{1}^{+2}\}) = \emptyset$. 
By the symmetry of $u_{1}$ and $u_{2}$, we also have $N_{G}(x_{1}, \{u_{2}^{+}; u_{2}^{+2}\}) = \emptyset$, 
and hence by the symmetry of $x_{1}$ and $x_{2}$, 
we have $N_{G}(x_{2}, \{u_{1}^{+}, u_{1}^{+2}, u_{2}^{+}, u_{2}^{+2}\}) = \emptyset$. 
In particular, 
since $\{x_{1}, x_{2}\} \subseteq A$, 
this implies that $\{u_{1}^{+}, u_{1}^{+2}, u_{2}^{+}, u_{2}^{+2}\} \subseteq A$. 
\qed

Since $u_{1} \in B$ and $\{x_{1}, x_{2}, u_{1}^{+}, u_{1}^{+2}, u_{2}^{+}, u_{2}^{+2}\} \subseteq A$ by Claim~\ref{cl:P4-2}, 
we have that 
$G[\{u_{1}, x_{1}, x_{2}, u_{1}^{+},$ $u_{1}^{+2}, u_{2}^{+}, u_{2}^{+2}\}]$ contains $W$ as a subgraph. 
Since 
$G$ is $W$-free, Claim~\ref{cl:P4-2} yields that 
$N_{G}(\{u_{1}^{+}, u_{1}^{+2}\};$ $\{u_{2}^{+}, u_{2}^{+2}\}) \neq \emptyset$. 
If $u_{1}^{+}u_{2}^{+} \in E(G)$, 
then 
$u_{1}^{+}u_{2}^{+}\ora{C}u_{1}x_{1}u_{2}\ola{C}u_{1}^{+}$ is a longer cycle than $C$, a contradiction. 
If $u_{1}^{+}u_{2}^{+2} \in E(G)$, 
then 
$u_{1}^{+}u_{2}^{+2}\ora{C}u_{1}x_{1}x_{2}u_{2}\ola{C}u_{1}^{+}$ is a longer cycle than $C$, a contradiction again. 
Similarly, we have $u_{1}^{+2}u_{2}^{+} \notin E(G)$, 
and hence we have $u_{1}^{+2}u_{2}^{+2} \in E(G)$. 
But, then $u_{1}^{+}u_{1}^{+2}u_{2}^{+2}u_{2}^{+}$ is an induced path of $G$, 
which contradicts the assumption that $G$ is $P_{4}$-free. 

This completes the proof of Theorem~\ref{P_{4}, W, long}.
\qed

\subsection{Proof of Theorem~\ref{K_{1,3}^{**}, Z_{1}}}
\label{Proof of K_{1, 3}^{**}, Z_{1}}

The proof of Theorem~\ref{K_{1,3}^{**}, Z_{1}} is actually divided into two parts according as 
the graph contains a triangle or not. 
To do that, we use the following.

\begin{Lem}[Olariu~\cite{O}]
\label{Olariu}
Let $G$ be a connected $Z_{1}$-free graph.
If $G$ contains a triangle, then $G$ is a complete multipartite graph.
\end{Lem}

\begin{thm}
\label{tri th}
A longest cycle of a $2$-connected $\{ K_{1, 3}^{**}, K_{3} \}$-free graph is a dominating cycle of the graph.
\end{thm}

Here we prove Theorem~\ref{K_{1,3}^{**}, Z_{1}} assuming Theorem~\ref{tri th}. 
We will show Theorem~\ref{tri th} in Subsections~\ref{longest cycles} and \ref{Proof of tri th}.

\medskip
\noindent
\textbf{Proof of Theorem~\ref{K_{1,3}^{**}, Z_{1}}.}~Let $G$ be a $2$-connected $\{K^{**}_{1,3},Z_{1}\}$-free graph.
If $G$ is $K_{3}$-free, then by Theorem~\ref{tri th}, 
$G$ has a longest cycle which is a dominating cycle.
Thus we may assume that $G$ contains a triangle.
Then by Lemma~\ref{Olariu}, $G$ is a complete multipartite graph. 
Let $\ora{C}$ be a longest of $G$.  
Suppose that there exists an edge $xy$ in $G - V(C)$, 
and let $u \in V(C)$. 
If some vertex $a$ in $\{x, y\}$ belongs to a different partite set from $u$ and $u^{+}$, 
then $uau^{+}\ora{C}u$ is a longer cycle than $C$, a contradiction. 
Thus we may assume that $u$ and $x$ belong to the same partite set, 
and $u^{+}$ and $y$ belong to the same partite set 
(note that $u$ and $u^{+}$ belong to different partite sets). 
Then $uyxu^{+}\ora{C}u$ is a longer cycle than $C$, a contradiction again. 
Thus $C$ is a dominating cycle of $G$. 
\qed

\subsection{Preparation for the proof of Theorem~\ref{tri th}}
\label{longest cycles}

In this subsection, 
we prepare lemmas which will be used in the proof of Theorem~\ref{tri th}.

Now let $G$ be a $2$-connected graph. 
Let $\ora{C}$ be a longest cycle of $G$, 
and let $H$ be a component of $G - V(C)$. 
(Note that $|V(C)| \ge 4$.)
Then by the maximality of $|V(C)|$, 
we can easily obtain the following lemma.

\begin{lem}
\label{maximality of C}
\begin{enumerate}[{\upshape(i)}]
\item $N_{G}(H; C) \cap N_{G}(H; C)^{+} = \emptyset$ (see the left of Figure~\ref{maximality-cycle}). 
\vspace{-5pt}
\item If $u$ and $v$ are distinct two vertices in $N_{G}(H; C)$, 
then 
$E(G) \cap \{u^{+}v^{+}, u^{-}v^{-} \} = \emptyset$ (see the center of Figure~\ref{maximality-cycle}). 
\vspace{-5pt}
\item If $u$ and $v$ are distinct two vertices in $N_{G}(H; C)$ 
such that $|N_{G}(u; H) \cup N_{G}(v; H)| \ge 2$, 
then $E(G) \cap \{u^{+2}v^{+}, u^{+}v^{+2}, u^{-2}v^{-}, u^{-}v^{-2} \} = \emptyset$. 
In particular, 
$u^{+2} \neq v$ and $u^{-2} \neq v^{-}$ (see the right of Figure~\ref{maximality-cycle}). 
\end{enumerate}
\end{lem}
\begin{figure}
\begin{center}
\unitlength 0.1in
\begin{picture}( 56.1500, 11.1500)(  11.5000,-24.9500)
%
\special{pn 12}%
\special{pa 2210 2000}%
\special{pa 2210 1560}%
\special{fp}%
%
\special{pn 12}%
\special{pa 2210 2000}%
\special{pa 2410 1600}%
\special{da 0.045}%
%
\special{pn 12}%
\special{ar 2206 2006 460 460  0.0000000 6.2831853}%
%
\special{pn 12}%
\special{ar 4006 2006 460 460  0.0000000 6.2831853}%
%
\special{pn 12}%
\special{pa 4006 1546}%
\special{pa 4006 2456}%
\special{fp}%
%
\special{pn 12}%
\special{ar 5806 2006 460 460  0.0000000 6.2831853}%
%
\special{pn 12}%
\special{sh 0}%
\special{ar 5806 2006 230 170  0.0000000 6.2831853}%
%
\special{pn 8}%
\special{sh 1.000}%
\special{ar 2206 1556 40 40  0.0000000 6.2831853}%
%
\special{pn 12}%
\special{sh 0}%
\special{ar 2206 2006 230 170  0.0000000 6.2831853}%
%
\special{pn 8}%
\special{sh 1.000}%
\special{ar 2406 1596 40 40  0.0000000 6.2831853}%
%
\special{pn 12}%
\special{sh 0}%
\special{ar 4006 2006 230 170  0.0000000 6.2831853}%
%
\special{pn 8}%
\special{sh 1.000}%
\special{ar 4006 1556 40 40  0.0000000 6.2831853}%
%
\special{pn 8}%
\special{sh 1.000}%
\special{ar 4006 2456 40 40  0.0000000 6.2831853}%
%
\special{pn 8}%
\special{sh 1.000}%
\special{ar 4206 1596 40 40  0.0000000 6.2831853}%
%
\special{pn 8}%
\special{sh 1.000}%
\special{ar 3806 2416 40 40  0.0000000 6.2831853}%
%
\special{pn 8}%
\special{sh 1.000}%
\special{ar 5806 1556 40 40  0.0000000 6.2831853}%
%
\special{pn 8}%
\special{sh 1.000}%
\special{ar 5806 2456 40 40  0.0000000 6.2831853}%
%
\special{pn 8}%
\special{sh 1.000}%
\special{ar 5606 2416 40 40  0.0000000 6.2831853}%
%
\special{pn 8}%
\special{sh 1.000}%
\special{ar 6006 1596 40 40  0.0000000 6.2831853}%
%
\special{pn 8}%
\special{sh 1.000}%
\special{ar 6166 1716 40 40  0.0000000 6.2831853}%
%
\special{pn 8}%
\special{sh 1.000}%
\special{ar 5806 1896 40 40  0.0000000 6.2831853}%
%
\special{pn 8}%
\special{sh 1.000}%
\special{ar 5806 2116 40 40  0.0000000 6.2831853}%
%
\special{pn 12}%
\special{pa 5806 2116}%
\special{pa 5806 2466}%
\special{fp}%
\special{pa 5806 1886}%
\special{pa 5806 1556}%
\special{fp}%
\put(56.7500,-19.9500){\makebox(0,0){{\scriptsize $H$}}}%
\put(58.0500,-14.6500){\makebox(0,0){{\scriptsize $u$}}}%
\put(61.0500,-15.0500){\makebox(0,0){{\scriptsize $u^{+}$}}}%
\put(63.1500,-16.5500){\makebox(0,0){{\scriptsize $u^{+2}$}}}%
\put(58.0500,-25.4500){\makebox(0,0){{\scriptsize $v$}}}%
\put(54.9500,-24.8500){\makebox(0,0){{\scriptsize $v^{+}$}}}%
\put(40.0500,-14.6500){\makebox(0,0){{\scriptsize $u$}}}%
\put(43.0500,-15.0500){\makebox(0,0){{\scriptsize $u^{+}$}}}%
\put(38.7500,-19.9500){\makebox(0,0){{\scriptsize $H$}}}%
\put(40.1000,-25.5000){\makebox(0,0){{\scriptsize $v$}}}%
\put(37.0000,-24.9000){\makebox(0,0){{\scriptsize $v^{+}$}}}%
\put(16.8500,-17.9500){\makebox(0,0){{\scriptsize $\ora{C}$}}}%
\put(20.7500,-19.9500){\makebox(0,0){{\scriptsize $H$}}}%
%
\special{pn 12}%
\special{pa 4206 1586}%
\special{pa 4218 1620}%
\special{pa 4232 1656}%
\special{pa 4244 1690}%
\special{pa 4256 1724}%
\special{pa 4268 1758}%
\special{pa 4278 1792}%
\special{pa 4288 1824}%
\special{pa 4296 1858}%
\special{pa 4302 1890}%
\special{pa 4308 1920}%
\special{pa 4312 1952}%
\special{pa 4314 1982}%
\special{pa 4312 2010}%
\special{pa 4310 2038}%
\special{pa 4304 2066}%
\special{pa 4298 2090}%
\special{pa 4286 2116}%
\special{pa 4274 2138}%
\special{pa 4258 2162}%
\special{pa 4240 2182}%
\special{pa 4220 2202}%
\special{pa 4196 2222}%
\special{pa 4172 2240}%
\special{pa 4146 2258}%
\special{pa 4118 2276}%
\special{pa 4090 2292}%
\special{pa 4060 2308}%
\special{pa 4028 2324}%
\special{pa 3994 2338}%
\special{pa 3962 2352}%
\special{pa 3928 2368}%
\special{pa 3892 2382}%
\special{pa 3858 2396}%
\special{pa 3822 2408}%
\special{pa 3806 2416}%
\special{sp 0.045}%
%
\special{pn 12}%
\special{pa 6166 1716}%
\special{pa 6162 1750}%
\special{pa 6160 1786}%
\special{pa 6158 1820}%
\special{pa 6154 1856}%
\special{pa 6150 1890}%
\special{pa 6146 1924}%
\special{pa 6140 1956}%
\special{pa 6134 1988}%
\special{pa 6126 2020}%
\special{pa 6118 2050}%
\special{pa 6108 2080}%
\special{pa 6098 2108}%
\special{pa 6084 2136}%
\special{pa 6070 2160}%
\special{pa 6054 2184}%
\special{pa 6036 2208}%
\special{pa 6016 2228}%
\special{pa 5994 2248}%
\special{pa 5970 2266}%
\special{pa 5944 2284}%
\special{pa 5918 2298}%
\special{pa 5890 2314}%
\special{pa 5860 2328}%
\special{pa 5830 2340}%
\special{pa 5798 2352}%
\special{pa 5766 2364}%
\special{pa 5734 2376}%
\special{pa 5700 2386}%
\special{pa 5668 2396}%
\special{pa 5634 2408}%
\special{pa 5606 2416}%
\special{sp 0.045}%
\put(17.4500,-15.0500){\makebox(0,0){{\bf {\small (i)}}}}%
\put(35.4500,-15.0500){\makebox(0,0){{\bf {\small (ii)}}}}%
\put(53.4500,-15.0500){\makebox(0,0){{\bf {\small (iii)}}}}%
\put(34.8500,-17.9500){\makebox(0,0){{\scriptsize $\ora{C}$}}}%
\put(52.8500,-17.9500){\makebox(0,0){{\scriptsize $\ora{C}$}}}%
\end{picture}%
\caption{Lemma~\ref{maximality of C}}
\label{maximality-cycle}
\end{center}
\end{figure}

Moreover, we give the following lemma concerning $\{K_{1, 3}^{**}, K_{3}\}$-free graphs.

\begin{lem}
\label{maximality of C and K_{1,3}^{**}Z_{1}-free}
Let $x \in V(H)$ and $u\in N_{G}(x;C)$.
If $G$ is $\{K_{1, 3}^{**}, K_{3}\}$-free and $|V(H)| \ge 2$,
then 
$E(G) \cap \{u^{+2}u^{-}, u^{+2}x\} \neq \emptyset$ and $E(G) \cap \{u^{-2}u^{+}, u^{-2}x\} \neq \emptyset$. 
\end{lem}

\noindent
\textbf{Proof of Lemma~\ref{maximality of C and K_{1,3}^{**}Z_{1}-free}.}~Since $|V(H)|\geq 2$, $N_{H}(x)\not=\emptyset $.
Let $x'\in N_{H}(x)$.
Then by Lemma~\ref{maximality of C} (i),
$E(G) \cap \{u^{+}x, u^{-}x, u^{+}x', u^{-}x'\} = \emptyset$.
By Lemma~\ref{maximality of C} (iii), $u^{+2}x' \notin E(G)$. 
Moreover, since $G$ is $K_{3}$-free,
$E(G) \cap \{uu^{+2}, u^{+}u^{-}, ux'\} = \emptyset$. 
Therefore, 
if 
$E(G) \cap \{u^{+2}u^{-}, u^{+2}x\} = \emptyset$, 
then $G[\{u, u^{+}, u^{+2}, u^{-}, x, x'\}]$ is isomorphic to $K_{1, 3}^{**}$, 
a contradiction. 
Thus $E(G) \cap \{u^{+2}u^{-}, u^{+2}x\} \neq \emptyset$. 
By the symmetry of $\ora{C}$ and $\ola{C}$, 
we have that $E(G) \cap \{u^{-2}u^{+}, u^{-2}x\} \neq \emptyset$. 
\qed

\subsection{Proof of Theorem~\ref{tri th}}
\label{Proof of tri th}

In this subsection, we prove Theorem~\ref{tri th}.

\medskip
\noindent
\textbf{Proof of Theorem~\ref{tri th}.}~Let $G$ be a $2$-connected $\{K_{1, 3}^{**}, K_{3}\}$-free graph, 
and 
we show that $G$ has a longest cycle which is a dominating cycle of $G$. 
By way of a contradiction, 
suppose that every longest cycle of $G$ is not a dominating cycle of $G$. 
For a cycle $C$ of $G$,
let $\mu(C) = \max \{|V(H)| : H$ is a component of $G - V(C)\}$.
Then $\mu(C) \ge 2$ for every longest cycle $C$ of $G$. 
For a cycle $C$ of $G$, we further define $\omega(C) = |\{H : H$ is a component of $G- V(C)$ such that $|V(H)| = \mu(C)\}|$.
Let $\ora{C}$ be a longest cycle of $G$. 
We choose $C$ so that 
\begin{enumerate}[]
\item{(C1)} $\mu(C)$ is as small as possible, and 
\vspace{-6pt} 
\item{(C2)} $\omega(C)$ is as small as possible, subject to (C1). 
\end{enumerate}

\noindent
Let $H$ be a component of $G- V(C)$ such that $|V(H)| = \mu(C) \ ( \ \ge 2)$. 
Since $G$ is $2$-connected, 
there exist distinct two vertices $u$ and $v$ in $C$ such that 
$N_{G}(u; H) \neq \emptyset$, $N_{G}(v; H) \neq \emptyset$, $|N_{G}(u; H) \cup N_{G}(v; H)| \ge 2$ 
and $N_{G}(H; u^{+} \ora{C} v^{-}) = \emptyset$. 
We choose the longest cycle $C$ of $G$, the component $H$ of $G-V(C)$ with $|V(H)| = \mu(C)$, the vertices $u$ and $v$ so that 
\begin{enumerate}[]
\item{(C3)} $|V(u^{+} \ora{C} v^{-})|$ is as large as possible, subject to (C1) and (C2). 
\end{enumerate}

By Lemma~\ref{maximality of C} (i) and (iii),
$|V(u \ora{C} v)|\geq 4$ and $|V(v \ora{C} u)|\geq 4$.
Since $N_{G}(H; u^{+} \ora{C} v^{-}) = \emptyset$,
it follows from Lemma~\ref{maximality of C and K_{1,3}^{**}Z_{1}-free} that 
$u^{+2}u^{-} \in E(G)$, 
and hence 
by Lemma~\ref{maximality of C} (ii) and (iii), $|V(u \ora{C} v)| \ge 6$.

\begin{claim}
\label{yu^{+3}}
$yu^{+3} \in E(G)$ for $y \in N_{G}(u^{+}; G-V(C))$.  
\end{claim}
\proof 
Suppose that $yu^{+3} \notin E(G)$ for some vertex $y\in N_{G}(u^{+}; G-V(C))$.
By the choice of $u$ and $v$, $y \notin V(H)$. 
Let $x \in N_{G}(u; H)$.
Since $G$ is $K_{3}$-free, $E(G) \cap \{yu, yu^{+2},uu^{+2},u^{+}u^{+3}\} = \emptyset$.
Since $N_{G}(H; u^{+} \ora{C} v^{-}) = \emptyset$, $yu^{+3} \notin E(G)$ and $G[\{x, y, u, u^{+}, u^{+2}, u^{+3}\}] \not \cong K_{1, 3}^{**}$, 
these imply that $uu^{+3} \in E(G)$ (see Figure~\ref{yu+3}).
However, $G[\{x, x', y, u, u^{+}, u^{+3}\}]$ is isomorphic to $K_{1, 3}^{**}$ where $x' \in N_{H}(x)$
because $N_{G}(H; u^{+} \ora{C} v^{-}) = \emptyset$ and $G$ is $K_{3}$-free, a contradiction. 
\qed

\begin{figure}
\begin{center}
\unitlength 0.1in
\begin{picture}( 27.4000, 13.6500)(  10.5000,-21.1000)
\put(14.6000,-12.8000){\makebox(0,0){{\small $\ora{C}$}}}%
%
\special{pn 7}%
\special{sh 1.000}%
\special{ar 2600 1400 50 50  0.0000000 6.2831853}%
%
\special{pn 7}%
\special{sh 1.000}%
\special{ar 3200 1400 50 50  0.0000000 6.2831853}%
%
\special{pn 7}%
\special{sh 1.000}%
\special{ar 1600 1400 50 50  0.0000000 6.2831853}%
%
\special{pn 11}%
\special{pa 1450 1400}%
\special{pa 1800 1400}%
\special{fp}%
%
\special{pn 22}%
\special{pa 1800 1400}%
\special{pa 2400 1400}%
\special{fp}%
%
\special{pn 11}%
\special{pa 2400 1400}%
\special{pa 2750 1400}%
\special{fp}%
%
\special{pn 11}%
\special{pa 3050 1400}%
\special{pa 3390 1400}%
\special{fp}%
%
\special{pn 11}%
\special{pa 3020 1400}%
\special{pa 2760 1400}%
\special{dt 0.045}%
%
\special{pn 22}%
\special{pa 2000 1400}%
\special{pa 2400 800}%
\special{fp}%
\put(32.5000,-12.9000){\makebox(0,0){{\small $v$}}}%
%
\special{pn 7}%
\special{pa 2200 1410}%
\special{pa 2400 800}%
\special{da 0.060}%
\special{pa 2400 800}%
\special{pa 2400 1400}%
\special{da 0.060}%
%
\special{pn 7}%
\special{pa 1800 1400}%
\special{pa 2400 800}%
\special{da 0.060}%
%
\special{pn 11}%
\special{pa 3210 1390}%
\special{pa 2880 1950}%
\special{fp}%
%
\special{pn 11}%
\special{sh 0}%
\special{ar 2510 1950 500 160  0.0000000 6.2831853}%
\put(25.3000,-19.6000){\makebox(0,0){{\small $H$}}}%
\put(22.3000,-19.8000){\makebox(0,0){{\small $x$}}}%
%
\special{pn 22}%
\special{pa 1800 1390}%
\special{pa 2130 1950}%
\special{fp}%
%
\special{pn 7}%
\special{pa 1800 1410}%
\special{pa 1816 1368}%
\special{pa 1832 1328}%
\special{pa 1846 1286}%
\special{pa 1862 1248}%
\special{pa 1878 1210}%
\special{pa 1894 1176}%
\special{pa 1908 1146}%
\special{pa 1924 1118}%
\special{pa 1940 1094}%
\special{pa 1956 1074}%
\special{pa 1972 1062}%
\special{pa 1986 1052}%
\special{pa 2002 1050}%
\special{pa 2018 1054}%
\special{pa 2034 1064}%
\special{pa 2050 1080}%
\special{pa 2066 1102}%
\special{pa 2082 1126}%
\special{pa 2098 1156}%
\special{pa 2114 1188}%
\special{pa 2130 1222}%
\special{pa 2146 1260}%
\special{pa 2162 1300}%
\special{pa 2178 1340}%
\special{pa 2194 1382}%
\special{pa 2200 1400}%
\special{sp 0.060}%
%
\special{pn 7}%
\special{pa 2130 1940}%
\special{pa 2000 1390}%
\special{da 0.060}%
\special{pa 2200 1390}%
\special{pa 2120 1940}%
\special{da 0.060}%
\special{pa 2120 1940}%
\special{pa 2400 1400}%
\special{da 0.060}%
%
\special{pn 7}%
\special{pa 2000 1400}%
\special{pa 2028 1438}%
\special{pa 2058 1474}%
\special{pa 2086 1510}%
\special{pa 2112 1544}%
\special{pa 2140 1576}%
\special{pa 2166 1606}%
\special{pa 2192 1634}%
\special{pa 2216 1658}%
\special{pa 2240 1676}%
\special{pa 2260 1692}%
\special{pa 2282 1702}%
\special{pa 2300 1708}%
\special{pa 2316 1706}%
\special{pa 2332 1700}%
\special{pa 2344 1686}%
\special{pa 2356 1666}%
\special{pa 2364 1642}%
\special{pa 2372 1612}%
\special{pa 2380 1578}%
\special{pa 2386 1542}%
\special{pa 2390 1504}%
\special{pa 2396 1462}%
\special{pa 2400 1422}%
\special{pa 2400 1410}%
\special{sp 0.060}%
%
\special{pn 11}%
\special{sh 0}%
\special{ar 2400 800 50 50  0.0000000 6.2831853}%
%
\special{pn 11}%
\special{sh 0}%
\special{ar 1800 1400 50 50  0.0000000 6.2831853}%
%
\special{pn 11}%
\special{sh 0}%
\special{ar 2000 1400 50 50  0.0000000 6.2831853}%
%
\special{pn 11}%
\special{sh 0}%
\special{ar 2200 1400 50 50  0.0000000 6.2831853}%
%
\special{pn 11}%
\special{sh 0}%
\special{ar 2400 1400 50 50  0.0000000 6.2831853}%
\put(17.4000,-15.0000){\makebox(0,0){{\small $u$}}}%
\put(20.2000,-14.9000){\makebox(0,0){{\small $u^{+}$}}}%
\put(22.5000,-15.0000){\makebox(0,0){{\small $u^{+2}$}}}%
\put(24.4000,-12.6000){\makebox(0,0){{\small $u^{+3}$}}}%
\put(25.0000,-8.3000){\makebox(0,0){{\small $y$}}}%
%
\special{pn 11}%
\special{sh 0}%
\special{ar 2120 1950 50 50  0.0000000 6.2831853}%
\end{picture}%
\caption{Claim~\ref{yu^{+3}}}
\label{yu+3}
\end{center}
\end{figure}

\begin{claim}
\label{N_{G}(u^{-2}; H) = emptyset}
$N_{G}(u^{-2}; H) = \emptyset$. 
\end{claim}
\proof 
Suppose that $N_{G}(u^{-2}; H) \neq \emptyset$, 
and let $x \in N_{G}(u; H)$. 
By Lemma~\ref{maximality of C} (iii), $N_{G}(u^{-2}; H) = \{x\}$. 
If there exists a vertex $y\in N_{G}(u^{+}; G-V(C))$, 
then by Claim~\ref{yu^{+3}}, $yu^{+3} \in E(G)$, 
and hence 
$u^{-2}xuu^{-}u^{+2}u^{+}yu^{+3}\ora{C}u^{-2}$ is a longer cycle than $C$ 
(note that by the choice of $u$ and $v$, $y \notin V(H)$), a contradiction.
Thus $N_{G}(u^{+}; G-V(C)) = \emptyset$. 
Then $D := u^{-2}xuu^{-}u^{+2}\ora{C}u^{-2}$ is a cycle in $G$ such that $V(D) = (V(C) \bss \{u^{+}\}) \cup \{x\}$. 
Since $N_{G}(u^{+}; G-V(C)) = \emptyset$, 
it follows that $u^{+}$ is  a component of $G-V(D)$, 
and $G-V(D)$ contains some components whose union is $H - \{x\}$, 
which contradicts the choice (C1) or (C2). 
\qed

By Lemma~\ref{maximality of C and K_{1,3}^{**}Z_{1}-free} 
and Claim~\ref{N_{G}(u^{-2}; H) = emptyset}, 
we have $u^{+}u^{-2} \in E(G)$, 
and hence by Lemma~\ref{maximality of C} (ii) and (iii), 
$|V(v \ora{C} u)| \ge 6$. 
Moreover, by the maximality of $|V(C)|$, 
we can easily see that the following holds.

\begin{claim}
\label{maximality of C (2)}
$E(G) \cap \{ u^{-3}v^{-}, u^{-4}v^{-} \} = \emptyset$. 
\end{claim}
\proof 
Let $x \in N_{G}(u; H)$ and $x' \in N_{G}(v; H)$ with $x \neq x'$, 
and let $\ora{P}$ be an $(x, x')$-path in $H$.
If $v^{-}$ is adjacent to a vertex $z\in \{u^{-3},u^{-4}\}$,
then $zv^{-} \ola{C} u^{+}u^{-2}$ $\ora{C} u x \ora{P} x' v \ora{C}z$
is a longer cycle than $C$, a contradiction.
Thus $E(G) \cap \{ u^{-3}v^{-}, u^{-4}v^{-} \} = \emptyset$. 
\qed

Let $w \in N_{G}(H; v \ora{C} u^{-})$. 
We choose $w$ so that $|V(w \ora{C}u)|$ is as small as possible. 
Note that by the choice of $w$, $N_{G}(H; w^{+} \ora{C} u^{-}) = \emptyset$ (possibly $w = v$). 
By Lemma~\ref{maximality of C} (i) and Claim~\ref{N_{G}(u^{-2}; H) = emptyset}, $w \notin \{ u^{-}, u^{-2} \}$. 
Since $u^{-2}u^{+} \in E(G)$, it follows from Lemma~\ref{maximality of C} (ii) that $w \neq u^{-3}$. 
Write $u^{-} \ola{C} w^{+} = u_{1}u_{2} \dots u_{s-1}u_{s}$ ($s \ge 3$). 
For an integer $k$ with $1 \le k \le s - 1$, 
we call $u_{1} \ola{C} u_{k} \ ( \ = u_{1}u_{2} \dots u_{k})$ an \textit{insertible path of} $C$ 
if there exist distinct $k$ vertices $v_{1}, v_{2}, \dots, v_{k}$ satisfying the following (see Figure~\ref{insertible}): 
\begin{figure}
\begin{center}
\unitlength 0.1in
\begin{picture}( 57.4000, 10.4100)( 19.3000,-21.6000)
%
\special{pn 8}%
\special{sh 1.000}%
\special{ar 3600 1600 50 50  0.0000000 6.2831853}%
%
\special{pn 8}%
\special{sh 1.000}%
\special{ar 3400 1600 50 50  0.0000000 6.2831853}%
%
\special{pn 8}%
\special{sh 1.000}%
\special{ar 3000 1600 50 50  0.0000000 6.2831853}%
%
\special{pn 8}%
\special{sh 1.000}%
\special{ar 3600 1600 50 50  0.0000000 6.2831853}%
%
\special{pn 8}%
\special{sh 1.000}%
\special{ar 3600 1600 50 50  0.0000000 6.2831853}%
%
\special{pn 8}%
\special{sh 1.000}%
\special{ar 3200 1600 50 50  0.0000000 6.2831853}%
%
\special{pn 8}%
\special{sh 1.000}%
\special{ar 3800 1600 50 50  0.0000000 6.2831853}%
%
\special{pn 8}%
\special{pa 3000 1600}%
\special{pa 3028 1628}%
\special{pa 3054 1654}%
\special{pa 3080 1678}%
\special{pa 3108 1702}%
\special{pa 3134 1724}%
\special{pa 3160 1746}%
\special{pa 3186 1762}%
\special{pa 3214 1778}%
\special{pa 3240 1790}%
\special{pa 3266 1796}%
\special{pa 3294 1800}%
\special{pa 3320 1800}%
\special{pa 3346 1794}%
\special{pa 3374 1784}%
\special{pa 3400 1772}%
\special{pa 3426 1754}%
\special{pa 3454 1736}%
\special{pa 3480 1714}%
\special{pa 3506 1692}%
\special{pa 3534 1666}%
\special{pa 3560 1642}%
\special{pa 3586 1614}%
\special{pa 3600 1600}%
\special{sp}%
%
\special{pn 8}%
\special{pa 3200 1600}%
\special{pa 3228 1628}%
\special{pa 3254 1654}%
\special{pa 3280 1678}%
\special{pa 3308 1702}%
\special{pa 3334 1724}%
\special{pa 3360 1746}%
\special{pa 3386 1762}%
\special{pa 3414 1778}%
\special{pa 3440 1790}%
\special{pa 3466 1796}%
\special{pa 3494 1800}%
\special{pa 3520 1800}%
\special{pa 3546 1794}%
\special{pa 3574 1784}%
\special{pa 3600 1772}%
\special{pa 3626 1754}%
\special{pa 3654 1736}%
\special{pa 3680 1714}%
\special{pa 3706 1692}%
\special{pa 3734 1666}%
\special{pa 3760 1642}%
\special{pa 3786 1614}%
\special{pa 3800 1600}%
\special{sp}%
%
\special{pn 8}%
\special{sh 1.000}%
\special{ar 2800 1600 50 50  0.0000000 6.2831853}%
%
\special{pn 8}%
\special{sh 1.000}%
\special{ar 2600 1600 50 50  0.0000000 6.2831853}%
%
\special{pn 8}%
\special{pa 2400 1600}%
\special{pa 4000 1600}%
\special{fp}%
%
\special{pn 8}%
\special{pa 3800 1600}%
\special{pa 3800 1600}%
\special{pa 3800 1600}%
\special{pa 3800 1600}%
\special{pa 3800 1600}%
\special{pa 3800 1600}%
\special{pa 3800 1600}%
\special{pa 3800 1600}%
\special{pa 3800 1600}%
\special{pa 3800 1600}%
\special{pa 3800 1600}%
\special{pa 3800 1600}%
\special{pa 3800 1600}%
\special{pa 3800 1600}%
\special{pa 3800 1600}%
\special{pa 3800 1600}%
\special{pa 3800 1600}%
\special{pa 3800 1600}%
\special{pa 3800 1600}%
\special{pa 3800 1600}%
\special{pa 3800 1600}%
\special{pa 3800 1600}%
\special{pa 3800 1600}%
\special{pa 3800 1600}%
\special{pa 3800 1600}%
\special{pa 3800 1600}%
\special{pa 3800 1600}%
\special{pa 3800 1600}%
\special{pa 3800 1600}%
\special{pa 3800 1600}%
\special{pa 3800 1600}%
\special{pa 3800 1600}%
\special{pa 3800 1600}%
\special{pa 3800 1600}%
\special{pa 3800 1600}%
\special{pa 3800 1600}%
\special{pa 3800 1600}%
\special{pa 3800 1600}%
\special{pa 3800 1600}%
\special{sp}%
%
\special{pn 8}%
\special{sh 1.000}%
\special{ar 4700 1600 50 50  0.0000000 6.2831853}%
%
\special{pn 8}%
\special{sh 1.000}%
\special{ar 6000 1600 50 50  0.0000000 6.2831853}%
%
\special{pn 8}%
\special{pa 3200 1600}%
\special{pa 3232 1590}%
\special{pa 3262 1578}%
\special{pa 3294 1568}%
\special{pa 3326 1558}%
\special{pa 3356 1548}%
\special{pa 3388 1538}%
\special{pa 3418 1528}%
\special{pa 3450 1518}%
\special{pa 3480 1510}%
\special{pa 3512 1502}%
\special{pa 3544 1494}%
\special{pa 3574 1486}%
\special{pa 3606 1480}%
\special{pa 3636 1474}%
\special{pa 3668 1470}%
\special{pa 3700 1466}%
\special{pa 3730 1464}%
\special{pa 3762 1462}%
\special{pa 3792 1460}%
\special{pa 3824 1460}%
\special{pa 3854 1462}%
\special{pa 3886 1464}%
\special{pa 3918 1468}%
\special{pa 3948 1472}%
\special{pa 3980 1478}%
\special{pa 4010 1484}%
\special{pa 4042 1490}%
\special{pa 4074 1498}%
\special{pa 4104 1506}%
\special{pa 4136 1514}%
\special{pa 4166 1522}%
\special{pa 4198 1532}%
\special{pa 4228 1542}%
\special{pa 4260 1552}%
\special{pa 4292 1562}%
\special{pa 4322 1574}%
\special{pa 4354 1584}%
\special{pa 4384 1594}%
\special{pa 4400 1600}%
\special{sp}%
%
\special{pn 8}%
\special{pa 3000 1600}%
\special{pa 3032 1590}%
\special{pa 3062 1578}%
\special{pa 3094 1566}%
\special{pa 3124 1556}%
\special{pa 3156 1544}%
\special{pa 3186 1534}%
\special{pa 3218 1524}%
\special{pa 3248 1514}%
\special{pa 3280 1504}%
\special{pa 3310 1494}%
\special{pa 3342 1484}%
\special{pa 3372 1474}%
\special{pa 3404 1466}%
\special{pa 3434 1458}%
\special{pa 3466 1450}%
\special{pa 3496 1442}%
\special{pa 3528 1436}%
\special{pa 3558 1428}%
\special{pa 3590 1422}%
\special{pa 3622 1418}%
\special{pa 3652 1414}%
\special{pa 3684 1410}%
\special{pa 3714 1406}%
\special{pa 3746 1404}%
\special{pa 3776 1402}%
\special{pa 3808 1400}%
\special{pa 3838 1400}%
\special{pa 3870 1400}%
\special{pa 3902 1402}%
\special{pa 3932 1404}%
\special{pa 3964 1406}%
\special{pa 3994 1410}%
\special{pa 4026 1414}%
\special{pa 4058 1418}%
\special{pa 4088 1424}%
\special{pa 4120 1430}%
\special{pa 4150 1436}%
\special{pa 4182 1442}%
\special{pa 4214 1450}%
\special{pa 4244 1458}%
\special{pa 4276 1466}%
\special{pa 4306 1474}%
\special{pa 4338 1482}%
\special{pa 4370 1492}%
\special{pa 4400 1500}%
\special{pa 4432 1510}%
\special{pa 4464 1520}%
\special{pa 4494 1530}%
\special{pa 4526 1540}%
\special{pa 4558 1552}%
\special{pa 4588 1562}%
\special{pa 4620 1572}%
\special{pa 4652 1584}%
\special{pa 4682 1594}%
\special{pa 4700 1600}%
\special{sp}%
%
\special{pn 8}%
\special{pa 3000 1600}%
\special{pa 3032 1588}%
\special{pa 3062 1574}%
\special{pa 3092 1560}%
\special{pa 3124 1546}%
\special{pa 3154 1534}%
\special{pa 3184 1520}%
\special{pa 3216 1508}%
\special{pa 3246 1494}%
\special{pa 3276 1482}%
\special{pa 3308 1470}%
\special{pa 3338 1456}%
\special{pa 3370 1444}%
\special{pa 3400 1432}%
\special{pa 3430 1420}%
\special{pa 3462 1408}%
\special{pa 3492 1398}%
\special{pa 3522 1386}%
\special{pa 3554 1376}%
\special{pa 3584 1366}%
\special{pa 3614 1356}%
\special{pa 3646 1346}%
\special{pa 3676 1336}%
\special{pa 3708 1328}%
\special{pa 3738 1318}%
\special{pa 3768 1310}%
\special{pa 3800 1302}%
\special{pa 3830 1296}%
\special{pa 3860 1288}%
\special{pa 3892 1282}%
\special{pa 3922 1276}%
\special{pa 3952 1272}%
\special{pa 3984 1266}%
\special{pa 4014 1262}%
\special{pa 4044 1258}%
\special{pa 4076 1256}%
\special{pa 4106 1254}%
\special{pa 4138 1252}%
\special{pa 4168 1250}%
\special{pa 4198 1250}%
\special{pa 4230 1250}%
\special{pa 4260 1252}%
\special{pa 4290 1254}%
\special{pa 4322 1256}%
\special{pa 4352 1258}%
\special{pa 4382 1262}%
\special{pa 4414 1266}%
\special{pa 4444 1270}%
\special{pa 4476 1276}%
\special{pa 4506 1282}%
\special{pa 4536 1288}%
\special{pa 4568 1294}%
\special{pa 4598 1302}%
\special{pa 4628 1310}%
\special{pa 4660 1318}%
\special{pa 4690 1326}%
\special{pa 4720 1334}%
\special{pa 4752 1344}%
\special{pa 4782 1354}%
\special{pa 4812 1364}%
\special{pa 4844 1374}%
\special{pa 4874 1386}%
\special{pa 4906 1396}%
\special{pa 4936 1408}%
\special{pa 4966 1418}%
\special{pa 4998 1430}%
\special{pa 5028 1442}%
\special{pa 5058 1456}%
\special{pa 5090 1468}%
\special{pa 5120 1480}%
\special{pa 5150 1492}%
\special{pa 5182 1506}%
\special{pa 5212 1518}%
\special{pa 5244 1532}%
\special{pa 5274 1546}%
\special{pa 5304 1558}%
\special{pa 5336 1572}%
\special{pa 5366 1586}%
\special{pa 5396 1598}%
\special{pa 5400 1600}%
\special{sp}%
%
\special{pn 8}%
\special{pa 2800 1600}%
\special{pa 2832 1588}%
\special{pa 2862 1574}%
\special{pa 2892 1560}%
\special{pa 2924 1548}%
\special{pa 2954 1534}%
\special{pa 2984 1522}%
\special{pa 3016 1508}%
\special{pa 3046 1496}%
\special{pa 3076 1484}%
\special{pa 3108 1470}%
\special{pa 3138 1458}%
\special{pa 3168 1446}%
\special{pa 3200 1434}%
\special{pa 3230 1422}%
\special{pa 3260 1410}%
\special{pa 3292 1398}%
\special{pa 3322 1386}%
\special{pa 3352 1374}%
\special{pa 3384 1364}%
\special{pa 3414 1352}%
\special{pa 3444 1342}%
\special{pa 3476 1332}%
\special{pa 3506 1322}%
\special{pa 3536 1312}%
\special{pa 3568 1302}%
\special{pa 3598 1292}%
\special{pa 3630 1284}%
\special{pa 3660 1276}%
\special{pa 3690 1266}%
\special{pa 3722 1258}%
\special{pa 3752 1252}%
\special{pa 3782 1244}%
\special{pa 3814 1238}%
\special{pa 3844 1232}%
\special{pa 3874 1226}%
\special{pa 3906 1220}%
\special{pa 3936 1214}%
\special{pa 3968 1210}%
\special{pa 3998 1206}%
\special{pa 4028 1202}%
\special{pa 4060 1198}%
\special{pa 4090 1196}%
\special{pa 4120 1194}%
\special{pa 4152 1192}%
\special{pa 4182 1190}%
\special{pa 4214 1190}%
\special{pa 4244 1190}%
\special{pa 4274 1190}%
\special{pa 4306 1192}%
\special{pa 4336 1194}%
\special{pa 4368 1196}%
\special{pa 4398 1198}%
\special{pa 4428 1200}%
\special{pa 4460 1204}%
\special{pa 4490 1208}%
\special{pa 4522 1212}%
\special{pa 4552 1218}%
\special{pa 4582 1222}%
\special{pa 4614 1228}%
\special{pa 4644 1234}%
\special{pa 4676 1240}%
\special{pa 4706 1248}%
\special{pa 4738 1254}%
\special{pa 4768 1262}%
\special{pa 4798 1270}%
\special{pa 4830 1278}%
\special{pa 4860 1288}%
\special{pa 4892 1296}%
\special{pa 4922 1306}%
\special{pa 4954 1314}%
\special{pa 4984 1324}%
\special{pa 5016 1334}%
\special{pa 5046 1344}%
\special{pa 5076 1356}%
\special{pa 5108 1366}%
\special{pa 5138 1376}%
\special{pa 5170 1388}%
\special{pa 5200 1400}%
\special{pa 5232 1410}%
\special{pa 5262 1422}%
\special{pa 5294 1434}%
\special{pa 5324 1446}%
\special{pa 5354 1458}%
\special{pa 5386 1470}%
\special{pa 5416 1484}%
\special{pa 5448 1496}%
\special{pa 5478 1508}%
\special{pa 5510 1520}%
\special{pa 5540 1534}%
\special{pa 5572 1546}%
\special{pa 5602 1560}%
\special{pa 5634 1572}%
\special{pa 5664 1586}%
\special{pa 5694 1598}%
\special{pa 5700 1600}%
\special{sp}%
%
\special{pn 8}%
\special{pa 2600 1600}%
\special{pa 2632 1588}%
\special{pa 2662 1574}%
\special{pa 2692 1560}%
\special{pa 2722 1548}%
\special{pa 2754 1534}%
\special{pa 2784 1520}%
\special{pa 2814 1508}%
\special{pa 2846 1494}%
\special{pa 2876 1482}%
\special{pa 2906 1468}%
\special{pa 2936 1456}%
\special{pa 2968 1444}%
\special{pa 2998 1430}%
\special{pa 3028 1418}%
\special{pa 3060 1406}%
\special{pa 3090 1394}%
\special{pa 3120 1382}%
\special{pa 3152 1370}%
\special{pa 3182 1358}%
\special{pa 3212 1346}%
\special{pa 3244 1334}%
\special{pa 3274 1324}%
\special{pa 3304 1312}%
\special{pa 3334 1302}%
\special{pa 3366 1292}%
\special{pa 3396 1280}%
\special{pa 3426 1270}%
\special{pa 3458 1260}%
\special{pa 3488 1252}%
\special{pa 3518 1242}%
\special{pa 3550 1232}%
\special{pa 3580 1224}%
\special{pa 3610 1216}%
\special{pa 3642 1208}%
\special{pa 3672 1200}%
\special{pa 3702 1192}%
\special{pa 3734 1184}%
\special{pa 3764 1178}%
\special{pa 3794 1170}%
\special{pa 3826 1164}%
\special{pa 3856 1158}%
\special{pa 3886 1154}%
\special{pa 3918 1148}%
\special{pa 3948 1144}%
\special{pa 3980 1140}%
\special{pa 4010 1136}%
\special{pa 4040 1132}%
\special{pa 4072 1128}%
\special{pa 4102 1126}%
\special{pa 4132 1124}%
\special{pa 4164 1122}%
\special{pa 4194 1120}%
\special{pa 4226 1120}%
\special{pa 4256 1120}%
\special{pa 4286 1120}%
\special{pa 4318 1120}%
\special{pa 4348 1122}%
\special{pa 4380 1122}%
\special{pa 4410 1124}%
\special{pa 4440 1126}%
\special{pa 4472 1130}%
\special{pa 4502 1132}%
\special{pa 4534 1136}%
\special{pa 4564 1140}%
\special{pa 4594 1144}%
\special{pa 4626 1148}%
\special{pa 4656 1154}%
\special{pa 4688 1160}%
\special{pa 4718 1164}%
\special{pa 4750 1170}%
\special{pa 4780 1178}%
\special{pa 4812 1184}%
\special{pa 4842 1192}%
\special{pa 4872 1198}%
\special{pa 4904 1206}%
\special{pa 4934 1214}%
\special{pa 4966 1222}%
\special{pa 4996 1230}%
\special{pa 5028 1240}%
\special{pa 5058 1248}%
\special{pa 5090 1258}%
\special{pa 5120 1266}%
\special{pa 5152 1276}%
\special{pa 5182 1286}%
\special{pa 5214 1296}%
\special{pa 5244 1306}%
\special{pa 5276 1318}%
\special{pa 5306 1328}%
\special{pa 5336 1338}%
\special{pa 5368 1350}%
\special{pa 5398 1362}%
\special{pa 5430 1372}%
\special{pa 5460 1384}%
\special{pa 5492 1396}%
\special{pa 5522 1408}%
\special{pa 5554 1420}%
\special{pa 5584 1432}%
\special{pa 5616 1444}%
\special{pa 5646 1456}%
\special{pa 5678 1468}%
\special{pa 5708 1480}%
\special{pa 5740 1494}%
\special{pa 5770 1506}%
\special{pa 5802 1518}%
\special{pa 5832 1532}%
\special{pa 5864 1544}%
\special{pa 5894 1556}%
\special{pa 5926 1570}%
\special{pa 5956 1582}%
\special{pa 5988 1596}%
\special{pa 6000 1600}%
\special{sp}%
%
\special{pn 8}%
\special{sh 1.000}%
\special{ar 6300 1600 50 50  0.0000000 6.2831853}%
%
\special{pn 8}%
\special{pa 4210 1600}%
\special{pa 4910 1600}%
\special{fp}%
%
\special{pn 8}%
\special{pa 5200 1600}%
\special{pa 6500 1600}%
\special{fp}%
%
\special{pn 8}%
\special{sh 1.000}%
\special{ar 7100 1600 50 50  0.0000000 6.2831853}%
%
\special{pn 8}%
\special{pa 6900 1600}%
\special{pa 7300 1600}%
\special{fp}%
\put(32.4000,-17.3000){\makebox(0,0){{\small $u_{1}$}}}%
\put(30.4000,-17.2000){\makebox(0,0){{\small $u_{2}$}}}%
\put(28.4000,-17.2000){\makebox(0,0){{\small $u_{3}$}}}%
\put(26.4000,-17.2000){\makebox(0,0){{\small $u_{4}$}}}%
\put(34.2000,-17.2000){\makebox(0,0){{\small $u$}}}%
\put(44.4000,-17.2000){\makebox(0,0){{\small $v_{1}$}}}%
\put(47.6000,-17.0000){\makebox(0,0){{\small $v_{1}^{+}$}}}%
\put(53.4000,-17.2000){\makebox(0,0){{\small $v_{2}$}}}%
\put(57.3000,-17.0000){\makebox(0,0){{\small $v_{2}^{+} = v_{3}$}}}%
\put(60.7000,-14.7000){\makebox(0,0){{\small $v_{3}^{+}$}}}%
%
\special{pn 8}%
\special{ar 5270 2000 1400 160  0.0000000 6.2831853}%
%
\special{pn 8}%
\special{pa 4100 2000}%
\special{pa 3400 1600}%
\special{fp}%
%
\special{pn 8}%
\special{pa 6420 2000}%
\special{pa 7120 1600}%
\special{fp}%
%
\special{pn 8}%
\special{sh 1.000}%
\special{ar 6440 2000 50 50  0.0000000 6.2831853}%
%
\special{pn 8}%
\special{sh 1.000}%
\special{ar 4090 2000 50 50  0.0000000 6.2831853}%
\put(52.9000,-20.0000){\makebox(0,0){{\small $H$}}}%
\put(71.0000,-17.3000){\makebox(0,0){{\small $v$}}}%
\put(23.7000,-14.9000){\makebox(0,0){{\small $\ora{C}$}}}%
%
\special{pn 8}%
\special{pa 4010 1600}%
\special{pa 4210 1600}%
\special{dt 0.045}%
%
\special{pn 8}%
\special{pa 4900 1600}%
\special{pa 5200 1600}%
\special{dt 0.045}%
%
\special{pn 8}%
\special{pa 6510 1600}%
\special{pa 6910 1600}%
\special{dt 0.045}%
%
\special{pn 13}%
\special{sh 0}%
\special{ar 4400 1600 50 50  0.0000000 6.2831853}%
%
\special{pn 13}%
\special{sh 0}%
\special{ar 5400 1600 50 50  0.0000000 6.2831853}%
%
\special{pn 13}%
\special{sh 0}%
\special{ar 5700 1600 50 50  0.0000000 6.2831853}%
\end{picture}%
\caption{The insertible path $u_{1}u_{2}u_{3}$ of $C$}
\label{insertible}
\end{center}
\end{figure}
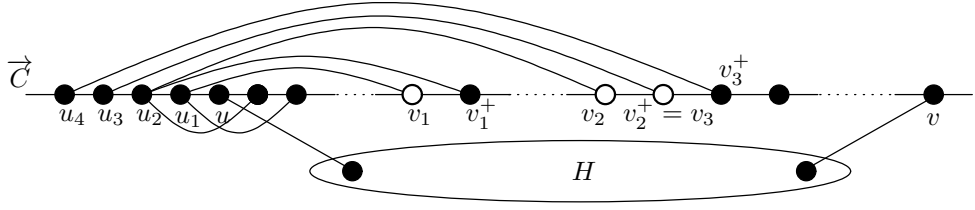

\begin{enumerate}[]
\item{(I1)} The vertices $v_{1}, v_{2}, \dots, v_{k}$ appear in this order along $\ora{C}$ and $v_{i} \in V(u^{+2} \ora{C} v^{-2})$ for each $i$ with $1 \le i \le k$. 
\vspace{-6pt}
\item{(I2)} $\{ u_{i}v_{i}, u_{i+1}v_{i}^{+} \} \subseteq E(G)$ for each $i$ with $1 \le i \le k$. 
\vspace{-6pt}
\item{(I3)} If $v_{k} \neq v^{-2}$, then $N_{G}(u_{i}; v_{i}^{+} \ora{C} v^{-2}) = \emptyset$ for each $i$ with $1 \le i \le k$. 
\end{enumerate}

\noindent
For an insertible path $u_{1} \ola{C} u_{k}$ of $C$, 
the vertices $v_{1}, v_{2}, \dots, v_{k}$ satisfying the conditions (I1)--(I3)
is called \textit{bridging vertices} of $u_{1} \ola{C} u_{k}$.

\begin{figure}
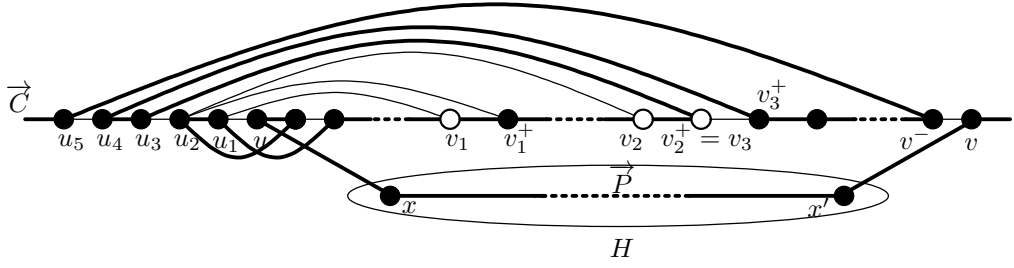

\begin{center}
\unitlength 0.1in
%
\caption{The case of $k-1 = 3$ and $l = 5$}
\label{l-odd}
\end{center}
\end{figure}

\begin{figure}
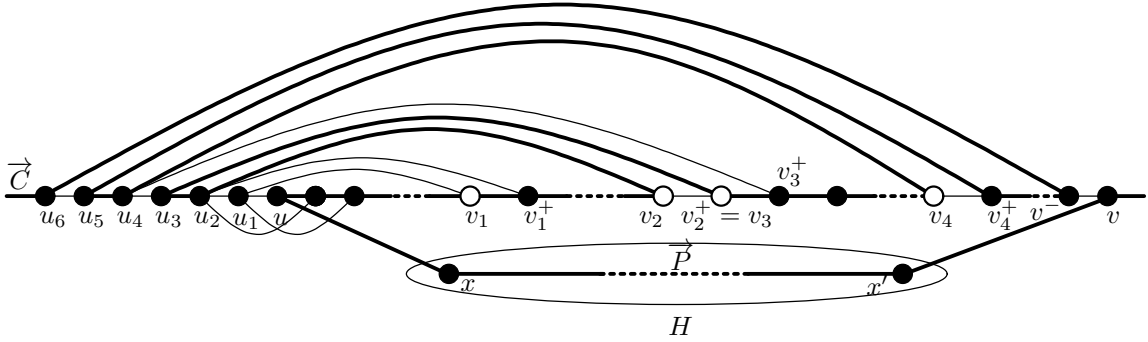

\begin{center}
\unitlength 0.1in
%
\caption{The case of $k-1 = 4$ and $l = 6$}
\label{l-even}
\end{center}
\end{figure}

\begin{claim}
\label{u_{l}v^{-}}
Let $k$ and $l$ be integers with $2 \le k \le s - 1$ and $1 \le l \le k+1$. 
If $u_{1} \ola{C} u_{k-1}$ is an insertible path of $C$, 
then $u_{l}v^{-} \notin E(G)$. 
In particular, if $v_{1}, \dots, v_{k-1}$ are bridging vertices of $u_{1} \ola{C} u_{k-1}$, then $v_{k-1} \neq v^{-2}$. 
\end{claim}
\proof 
By Lemma~\ref{maximality of C} (ii), (iii) and Claim~\ref{maximality of C (2)}, 
we may assume that 
$k \ge 4$ and $l \ge 5$.
Let $v_{1}, \dots, v_{k-1}$ be bridging vertices of $u_{1} \ola{C} u_{k-1}$.
Suppose that $u_{l}v^{-}\in E(G)$.
Let $x \in N_{G}(u; H)$ and $x' \in N_{G}(v; H)$ with $x \neq x'$, 
and let $\ora{P}$ be an $(x, x')$-path in $H$. 
If $l$ is odd, 
then by the condition (I2), 
$D := u_{l}v^{-} \ola{C} v_{l-2}^{+}u_{l-1}u_{l-2}v_{l-2} \ola{C} v_{l-4}^{+}u_{l-3}$ 
$u_{l-4}v_{l-4} \ola{C} \dots v_{3} \ola{C} u^{+2} u_{1}u_{2}u^{+}ux \ora{P}x' v \ora{C} u_{l}$
is a cycle in $G$ such that $V(D) = V(C) \cup V(P)$, which contradicts the maximality of $|V(C)|$ 
(see Figure~\ref{l-odd}). 
If $l$ is even, then 
$D := u_{l}v^{-} \ola{C} v_{l-2}^{+}u_{l-1}u_{l-2}$ $v_{l-2} \ola{C} v_{l-4}^{+}u_{l-3}u_{l-4}v_{l-4} \ola{C} \dots v_{2} \ola{C} ux \ora{P}x' v \ora{C} u_{l}$
is a cycle in $G$ such that $V(D) = (V(C) \bss \{ u_{1} \} ) \cup V(P)$, 
which contradicts the maximality of $|V(C)|$ again (see Figure~\ref{l-even}).  
\qed


\begin{claim}
\label{u_{k+2}u_{k}, u_{k+2}u_{k-1}}
Let $k$ be an integer with $2 \le k \le s - 1$. 
If $u_{1} \ola{C} u_{k-1}$ be an insertible path of $C$, 
then $u_{k+1}u_{k-2} \notin E(G)$, 
where $u_{0} = u$. 
\end{claim}
\proof 
Let $v_{1}, \dots, v_{k-1}$ be bridging vertices of $u_{1} \ola{C} u_{k-1}$. 
Suppose that $u_{k+1}u_{k - 2} \in E(G)$, 
and let $\ora{D} = u_{k+1}u_{k-2} \ora{C} v_{k-1}u_{k-1}u_{k}v_{k-1}^{+} \ora{C}u_{k+1}$. 
Then $D$ is a cycle in $G$ such that $V(D)  = V(C)$, 
and hence $\mu(D) = \mu(C)$ and $\omega(D) = \omega(C)$, 
in particular, $H$ is also a component of $G-V(D)$. 
Since $N_{G}(H; w^{+} \ora{C} u^{-} \cup u^{+} \ora{C} v^{-}) = \emptyset$, 
it follows from the definition of $D$ 
that $u$ and $v$ are distinct two vertices in $D$ such that 
$N_{G}(u; H) \neq \emptyset$, $N_{G}(v; H) \neq \emptyset$, $|N_{G}(u; H) \cup N_{G}(v; H)| \ge 2$ 
and $N_{G}(H; u^{+} \ora{D} v^{-}) = \emptyset$. 
Since $|V(u \ora{D} v)| > |V(u \ora{C} v)|$, this contradicts the choice (C3). 
\qed

\begin{figure}
\begin{center}
\unitlength 0.1in
\begin{picture}( 31.3000, 10.8000)( 14.0000,-19.5000)
\put(18.8000,-12.8000){\makebox(0,0){{\small $\ora{C}$}}}%
%
\special{pn 8}%
\special{sh 1.000}%
\special{ar 2600 1400 50 50  0.0000000 6.2831853}%
%
\special{pn 8}%
\special{sh 1.000}%
\special{ar 3200 1400 50 50  0.0000000 6.2831853}%
%
\special{pn 8}%
\special{sh 1.000}%
\special{ar 4000 1400 50 50  0.0000000 6.2831853}%
%
\special{pn 12}%
\special{pa 3900 1400}%
\special{pa 4200 1400}%
\special{fp}%
%
\special{pn 12}%
\special{pa 2790 1400}%
\special{pa 3010 1400}%
\special{dt 0.045}%
%
\special{pn 12}%
\special{pa 3770 1400}%
\special{pa 3870 1400}%
\special{dt 0.045}%
%
\special{pn 12}%
\special{pa 3670 1800}%
\special{pa 4000 1400}%
\special{fp}%
%
\special{pn 20}%
\special{pa 2200 1400}%
\special{pa 2230 1380}%
\special{pa 2258 1358}%
\special{pa 2286 1336}%
\special{pa 2314 1316}%
\special{pa 2344 1296}%
\special{pa 2372 1276}%
\special{pa 2400 1256}%
\special{pa 2430 1238}%
\special{pa 2458 1220}%
\special{pa 2486 1202}%
\special{pa 2516 1186}%
\special{pa 2544 1172}%
\special{pa 2572 1158}%
\special{pa 2602 1144}%
\special{pa 2630 1134}%
\special{pa 2658 1124}%
\special{pa 2688 1116}%
\special{pa 2716 1110}%
\special{pa 2744 1104}%
\special{pa 2772 1102}%
\special{pa 2802 1100}%
\special{pa 2830 1102}%
\special{pa 2858 1104}%
\special{pa 2888 1110}%
\special{pa 2916 1116}%
\special{pa 2944 1124}%
\special{pa 2974 1134}%
\special{pa 3002 1146}%
\special{pa 3030 1158}%
\special{pa 3060 1172}%
\special{pa 3088 1188}%
\special{pa 3116 1204}%
\special{pa 3146 1220}%
\special{pa 3174 1238}%
\special{pa 3202 1258}%
\special{pa 3230 1276}%
\special{pa 3260 1296}%
\special{pa 3288 1318}%
\special{pa 3316 1338}%
\special{pa 3346 1360}%
\special{pa 3374 1380}%
\special{pa 3400 1400}%
\special{sp}%
%
\special{pn 8}%
\special{pa 2200 1400}%
\special{pa 2228 1376}%
\special{pa 2256 1352}%
\special{pa 2284 1330}%
\special{pa 2312 1306}%
\special{pa 2340 1282}%
\special{pa 2368 1260}%
\special{pa 2394 1238}%
\special{pa 2422 1216}%
\special{pa 2450 1196}%
\special{pa 2478 1174}%
\special{pa 2506 1156}%
\special{pa 2534 1136}%
\special{pa 2562 1118}%
\special{pa 2590 1102}%
\special{pa 2618 1086}%
\special{pa 2646 1070}%
\special{pa 2672 1058}%
\special{pa 2700 1044}%
\special{pa 2728 1034}%
\special{pa 2756 1024}%
\special{pa 2784 1016}%
\special{pa 2812 1010}%
\special{pa 2840 1004}%
\special{pa 2868 1002}%
\special{pa 2896 1000}%
\special{pa 2922 1002}%
\special{pa 2950 1004}%
\special{pa 2978 1008}%
\special{pa 3006 1014}%
\special{pa 3034 1020}%
\special{pa 3062 1030}%
\special{pa 3090 1040}%
\special{pa 3118 1052}%
\special{pa 3146 1066}%
\special{pa 3172 1080}%
\special{pa 3200 1096}%
\special{pa 3228 1112}%
\special{pa 3256 1130}%
\special{pa 3284 1148}%
\special{pa 3312 1168}%
\special{pa 3340 1188}%
\special{pa 3368 1208}%
\special{pa 3396 1230}%
\special{pa 3422 1252}%
\special{pa 3450 1274}%
\special{pa 3478 1298}%
\special{pa 3506 1320}%
\special{pa 3534 1344}%
\special{pa 3562 1368}%
\special{pa 3590 1392}%
\special{pa 3600 1400}%
\special{sp 0.050}%
%
\special{pn 8}%
\special{pa 2000 1400}%
\special{pa 2028 1374}%
\special{pa 2054 1348}%
\special{pa 2080 1322}%
\special{pa 2108 1296}%
\special{pa 2134 1270}%
\special{pa 2160 1244}%
\special{pa 2188 1218}%
\special{pa 2214 1194}%
\special{pa 2240 1170}%
\special{pa 2268 1146}%
\special{pa 2294 1122}%
\special{pa 2320 1100}%
\special{pa 2348 1078}%
\special{pa 2374 1058}%
\special{pa 2400 1036}%
\special{pa 2428 1018}%
\special{pa 2454 998}%
\special{pa 2480 982}%
\special{pa 2508 966}%
\special{pa 2534 950}%
\special{pa 2560 936}%
\special{pa 2588 922}%
\special{pa 2614 912}%
\special{pa 2640 900}%
\special{pa 2668 892}%
\special{pa 2694 884}%
\special{pa 2720 878}%
\special{pa 2748 874}%
\special{pa 2774 872}%
\special{pa 2800 870}%
\special{pa 2828 872}%
\special{pa 2854 874}%
\special{pa 2882 878}%
\special{pa 2908 884}%
\special{pa 2934 890}%
\special{pa 2962 900}%
\special{pa 2988 910}%
\special{pa 3014 920}%
\special{pa 3042 934}%
\special{pa 3068 948}%
\special{pa 3096 962}%
\special{pa 3122 978}%
\special{pa 3148 996}%
\special{pa 3176 1014}%
\special{pa 3202 1034}%
\special{pa 3228 1054}%
\special{pa 3256 1074}%
\special{pa 3282 1096}%
\special{pa 3310 1118}%
\special{pa 3336 1142}%
\special{pa 3362 1164}%
\special{pa 3390 1188}%
\special{pa 3416 1214}%
\special{pa 3444 1238}%
\special{pa 3470 1264}%
\special{pa 3496 1290}%
\special{pa 3524 1316}%
\special{pa 3550 1342}%
\special{pa 3578 1368}%
\special{pa 3604 1394}%
\special{pa 3610 1400}%
\special{sp 0.050}%
%
\special{pn 8}%
\special{pa 2000 1400}%
\special{pa 2024 1434}%
\special{pa 2046 1468}%
\special{pa 2068 1498}%
\special{pa 2092 1526}%
\special{pa 2114 1552}%
\special{pa 2136 1572}%
\special{pa 2158 1588}%
\special{pa 2182 1598}%
\special{pa 2204 1600}%
\special{pa 2226 1596}%
\special{pa 2250 1584}%
\special{pa 2272 1566}%
\special{pa 2294 1544}%
\special{pa 2318 1518}%
\special{pa 2340 1488}%
\special{pa 2362 1456}%
\special{pa 2386 1424}%
\special{pa 2400 1400}%
\special{sp 0.050}%
%
\special{pn 12}%
\special{sh 0}%
\special{ar 3200 1800 570 150  0.0000000 6.2831853}%
%
\special{pn 20}%
\special{pa 2400 1400}%
\special{pa 2750 1800}%
\special{fp}%
%
\special{pn 8}%
\special{pa 2750 1800}%
\special{pa 2000 1400}%
\special{da 0.050}%
%
\special{pn 8}%
\special{pa 2760 1800}%
\special{pa 2200 1400}%
\special{da 0.050}%
\special{pa 3600 1400}%
\special{pa 2750 1800}%
\special{da 0.050}%
\special{pa 2750 1800}%
\special{pa 3410 1400}%
\special{da 0.050}%
%
\special{pn 12}%
\special{sh 0}%
\special{ar 2750 1800 50 50  0.0000000 6.2831853}%
\put(32.0000,-18.1000){\makebox(0,0){{\small $H$}}}%
\put(28.6000,-18.3000){\makebox(0,0){{\small $x$}}}%
\put(19.5000,-14.9000){\makebox(0,0){{\small $u_{2}$}}}%
\put(21.6000,-15.0000){\makebox(0,0){{\small $u_{1}$}}}%
\put(24.0000,-15.0000){\makebox(0,0){{\small $u$}}}%
\put(33.5000,-14.9000){\makebox(0,0){{\small $v_{1}$}}}%
\put(36.9000,-12.6000){\makebox(0,0){{\small $v_{1}^{+}$}}}%
\put(40.1000,-12.8000){\makebox(0,0){{\small $v$}}}%
%
\special{pn 20}%
\special{pa 2000 1400}%
\special{pa 2400 1400}%
\special{fp}%
%
\special{pn 20}%
\special{pa 3400 1400}%
\special{pa 3600 1400}%
\special{fp}%
%
\special{pn 12}%
\special{pa 3600 1400}%
\special{pa 3760 1400}%
\special{fp}%
%
\special{pn 12}%
\special{pa 3400 1400}%
\special{pa 3040 1400}%
\special{fp}%
%
\special{pn 12}%
\special{pa 2410 1400}%
\special{pa 2770 1400}%
\special{fp}%
%
\special{pn 12}%
\special{pa 2000 1400}%
\special{pa 1800 1400}%
\special{fp}%
%
\special{pn 12}%
\special{sh 0}%
\special{ar 2200 1400 50 50  0.0000000 6.2831853}%
%
\special{pn 12}%
\special{sh 0}%
\special{ar 3600 1400 50 50  0.0000000 6.2831853}%
%
\special{pn 8}%
\special{pa 2000 1400}%
\special{pa 2028 1376}%
\special{pa 2054 1352}%
\special{pa 2082 1328}%
\special{pa 2108 1304}%
\special{pa 2136 1280}%
\special{pa 2162 1258}%
\special{pa 2190 1234}%
\special{pa 2216 1212}%
\special{pa 2244 1192}%
\special{pa 2272 1170}%
\special{pa 2298 1150}%
\special{pa 2326 1132}%
\special{pa 2352 1114}%
\special{pa 2380 1096}%
\special{pa 2408 1080}%
\special{pa 2434 1066}%
\special{pa 2462 1052}%
\special{pa 2490 1040}%
\special{pa 2516 1030}%
\special{pa 2544 1020}%
\special{pa 2572 1012}%
\special{pa 2600 1006}%
\special{pa 2626 1002}%
\special{pa 2654 1000}%
\special{pa 2682 1000}%
\special{pa 2710 1002}%
\special{pa 2738 1004}%
\special{pa 2766 1008}%
\special{pa 2794 1016}%
\special{pa 2822 1024}%
\special{pa 2850 1032}%
\special{pa 2878 1044}%
\special{pa 2906 1056}%
\special{pa 2934 1068}%
\special{pa 2962 1082}%
\special{pa 2990 1098}%
\special{pa 3020 1114}%
\special{pa 3048 1132}%
\special{pa 3076 1150}%
\special{pa 3104 1170}%
\special{pa 3132 1190}%
\special{pa 3160 1210}%
\special{pa 3190 1232}%
\special{pa 3218 1254}%
\special{pa 3246 1276}%
\special{pa 3274 1298}%
\special{pa 3304 1320}%
\special{pa 3332 1344}%
\special{pa 3360 1368}%
\special{pa 3388 1390}%
\special{pa 3400 1400}%
\special{sp 0.050}%
%
\special{pn 8}%
\special{pa 2400 1400}%
\special{pa 2430 1416}%
\special{pa 2462 1432}%
\special{pa 2492 1446}%
\special{pa 2522 1462}%
\special{pa 2552 1476}%
\special{pa 2582 1490}%
\special{pa 2612 1502}%
\special{pa 2644 1516}%
\special{pa 2674 1526}%
\special{pa 2704 1536}%
\special{pa 2734 1546}%
\special{pa 2764 1554}%
\special{pa 2794 1560}%
\special{pa 2824 1566}%
\special{pa 2856 1568}%
\special{pa 2886 1570}%
\special{pa 2916 1570}%
\special{pa 2946 1568}%
\special{pa 2976 1564}%
\special{pa 3006 1558}%
\special{pa 3036 1552}%
\special{pa 3066 1544}%
\special{pa 3096 1534}%
\special{pa 3126 1524}%
\special{pa 3156 1512}%
\special{pa 3186 1498}%
\special{pa 3216 1484}%
\special{pa 3246 1470}%
\special{pa 3276 1456}%
\special{pa 3308 1440}%
\special{pa 3338 1424}%
\special{pa 3368 1408}%
\special{pa 3398 1392}%
\special{pa 3400 1390}%
\special{sp 0.050}%
%
\special{pn 12}%
\special{sh 0}%
\special{ar 2000 1400 50 50  0.0000000 6.2831853}%
%
\special{pn 12}%
\special{sh 0}%
\special{ar 2400 1400 50 50  0.0000000 6.2831853}%
%
\special{pn 12}%
\special{sh 0}%
\special{ar 3400 1400 50 50  0.0000000 6.2831853}%
\end{picture}%
\caption{$G[\{x, u, u_{1}, u_{2}, v_{1}, v_{1}^{+}\}]$}
\label{u1Cu1}
\end{center}
\end{figure}
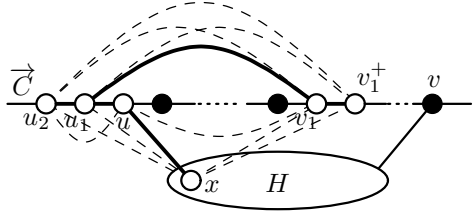

\begin{claim}
\label{u_{1}Cu_{k} is insertible}
For each $k$ with $1 \le k \le s - 1$, 
$u_{1} \ola{C} u_{k}$ is an insertible path of $C$.    
\end{claim}
\proof 
We first show that $u_{1} \ola{C} u_{1}  \ ( \ = u_{1})$ is an insertible path of $C$.
Since $(u^{-}u^{+2} = \ ) \ u_{1}u^{+2} \in E(G)$ 
and $|V(u \ora{C} v)| \ge 6$, there exists a vertex $v_{1}$ in $N_{G}(u_{1}; u^{+2} \ora {C} v^{-2})$. 
We choose $v_{1}$ so that $|V(v_{1} \ora{C} v^{-2})|$ is as small as possible. 
By Lemma~\ref{maximality of C} (ii), (iii) and the choice of $v_{1}$, 
we have
$v_{1} \in V(u^{+2} \ora {C} v^{-3})$ 
and 
$N_{G}(u_{1}; v_{1}^{+} \ora{C} v^{-}) = \emptyset$.
Suppose that $u_{2}v_{1}^{+} \notin E(G)$. 
Let $x \in N_{G}(u; H)$. 
Since $N_{G}(H; w^{+}\ora{C}u^{-} \cup u^{+}\ora{C}v^{-}) = \emptyset$ and $G$ is $K_{3}$-free, 
we have $E(G) \cap \{ xu_{2}, xu_{1}, xv_{1}, xv_{1}^{+}, u_{2}u, u_{2}v_{1},$ $uv_{1} \} = \emptyset$. 
Hence $G[\{x, u, u_{1}, u_{2}, v_{1}, v_{1}^{+}\}] \not \cong K_{1,3}^{**}$ yields that $uv_{1}^{+} \in E(G)$ 
(see Figure~\ref{u1Cu1}). 

Let $x' \in N_{H}(x)$. 
Recall that $v_{1} \in V(u^{+2} \ora {C} v^{-3})$.
By again the fact that $N_{G}(H; w^{+} \ora{C}u^{-} \cup u^{+} \ora{C} v^{-}) = \emptyset$ and $G$ is $K_{3}$-free,
we have that 
$E(G) \cap \{ xu_{1}, xv_{1}^{+}, xv_{1}^{+2}, x'u_{1}, x'u, x'v_{1}^{+},x'v_{1}^{+2},$ $uv_{1}^{+2} \} = \emptyset$. 
This together with the fact that $N_{G}(u_{1}) \cap  V(v_{1}^{+} \ora{C} v^{-}) = \emptyset$
implies that $G[\{x, x', u, u_{1},$ $v_{1}^{+}, v_{1}^{+2}\}] \cong K_{1, 3}^{**}$, a contradiction.
Thus $u_{2}v_{1}^{+} \in E(G)$, and hence $u_{1} \ola{C} u_{1}$ is an insatiable path of $C$. 

We next show that for $k$ with $2 \le k \le s - 1$, 
$u_{1} \ola{C} u_{k}$ is an insertible path of $C$.    
Suppose that there exists an integer $k$ with $2 \le k \le s - 1$
such that $u_{1} \ola{C} u_{k}$ is not an insertible path of $C$. 
We choose $k$ so that $k$ is as small as possible.
Then $u_{1} \ola{C} u_{k-1}$ is an insertible path of $C$. 
Since $u_{1} \ola{C} u_{k-1}$ is an insertible path of $C$, 
there exist bridging vertices $v_{1}, \dots, v_{k-1}$ of $u_{1} \ola{C} u_{k-1}$. 
Note that by Claim~\ref{u_{l}v^{-}}, $v_{k-1} \in V(u^{+2} \ora{C} v^{-3})$, 
and hence by the condition (I3) and again Claim~\ref{u_{l}v^{-}},  
we have 

\noindent
\begin{align}
\label{u_{k-1}, u_{k-2}} 
N_{G}( u_{i} ; v_{i}^{+} \ora{C} v^{-}) = \emptyset \textup{ for } 1 \le i \le k - 1.
\end{align}

\noindent
Since 
$v_{k-1} \in V(u^{+2} \ora{C} v^{-3})$, 
and since 
$u_{k}v_{k-1}^{+} \in E(G)$ by the condition (I2), 
there exists a vertex $v_{k}$ in $N_{G}(u_{k}; v_{k-1}^{+} \ora {C} v^{-2})$. 
We choose $v_{k}$ so that $|V(v_{k} \ora{C} v^{-2})|$ is as small as possible. 
Then the choice of $v_{k}$ implies that 
$N_{G}(u_{k}; v_{k}^{+} \ora{C} v^{-2}) = \emptyset$ if $v_{k} \neq v^{-2}$. 
Therefore, 
since $u_{1} \ola{C} u_{k-1}$ is an insertible path of $C$ 
and 
$u_{1} \ola{C} u_{k}$ is not an insertible path of $C$, 
we have 

\noindent
\begin{align}
\label{u_{k+1}, u_{k}} 
u_{k+1}v_{k}^{+} \notin E(G). 
\end{align} 

\noindent
Since $u_{1} \ola{C} u_{k-1}$ is an insertible path of $C$, 
it follows from Claim~\ref{u_{k+2}u_{k}, u_{k+2}u_{k-1}} that 
\begin{align}
\label{u_{k+1}, u_{k-2}}
u_{k+1} u_{k-2} \notin E(G), 
\end{align}
where $u_{0} = u$.
Since $G$ is $K_{3}$-free, 
we also have that

\noindent
\begin{align}
\label{K_{3}-free} 
E(G) \cap \{ u_{k+1}u_{k-1}, u_{k+1}v_{k}, u_{k}u_{k-2}, u_{k}v_{k}^{+} \} = \emptyset. 
\end{align}

If $k \ge 3$, then by combining (\ref{u_{k-1}, u_{k-2}})--(\ref{K_{3}-free}),  
we have 
$G[\{ u_{k+1}, u_{k}, u_{k-1}, u_{k-2}, v_{k}, v_{k}^{+} \}] \cong K_{1, 3}^{**}$, a contradiction (see the left of Figure~\ref{u1Cuk}).
Thus $k=2$.
Then by (\ref{u_{k-1}, u_{k-2}})--(\ref{K_{3}-free}), 
and 
since $G[\{ u_{3}, u_{2}, u_{1}, u, v_{2},$ $v_{2}^{+} \}] \not \cong K_{1, 3}^{**}$, 
we have $E(G)\cap \{uv_{2},uv_{2}^{+}\}\not=\emptyset $ (see the right of Figure~\ref{u1Cuk}).

\begin{figure}
\begin{center}
\unitlength 0.1in
%
\caption{Claim~\ref{u_{1}Cu_{k} is insertible}}
\label{u1Cuk}
\end{center}
\end{figure}

Let $x \in N_{G}(u; H)$ and $x' \in N_{H}(x)$.
If $uv_{2}^{+} \in E(G)$,
then since $N_{G}(H; w^{+}\ora{C}u^{-} \cup u^{+} \ora{C} v^{-}) = \emptyset$, 
it follows from (\ref{u_{k-1}, u_{k-2}}) and (\ref{K_{3}-free}) that $G[\{ x, x', u, u_{1}, u_{2}, v_{2}^{+} \}] \cong K_{1, 3}^{**}$, a contradiction. 
Thus $uv_{2}^{+} \notin E(G)$, and hence $uv_{2} \in E(G)$.
Then since $N_{G}(H; w^{+}\ora{C}u^{-} \cup u^{+} \ora{C} v^{-}) = \emptyset$, 
it follows from (\ref{u_{k-1}, u_{k-2}}) and (\ref{K_{3}-free}) that $G[\{ x, x', u, u_{1}, v_{2}, v_{2}^{+} \}] \cong K_{1, 3}^{**}$, a contradiction. 
\qed

By Claim~\ref{u_{1}Cu_{k} is insertible}, 
$u_{1} \ola{C} u_{s-1}$ is an insertible path of $C$. 
Let $v_{1}, \dots, v_{s-1}$ be bridging vertices of $u_{1} \ola{C} u_{s-1}$. 
Note that $v_{i} \in V(u^{+3} \ora{C}v^{-2})$ for each $i$ with $2 \le i \le s-1$. 
Let $x \in N_{G}(u; H)$ and $x' \in N_{G}(w; H)$ (if possible, choose $x'\not=x$), 
and let $\ora{P}$ be an $(x', x)$-path in $H$. 
Recall that $\{u_{1}u^{+2}, u_{2}u^{+} \} \subseteq E(G)$. 
If $s$ is even, then  
$D := wx' \ora{P} xu u^{+} u_{2} u_{1} u^{+2} \ora{C} v_{3} u_{3} u_{4}$ $v_{3}^{+} \ora{C} v_{5} u_{5}u_{6} v_{5}^{+} \ora{C} \dots $ $v_{s-1}u_{s-1}u_{s}v_{s-1}^{+} \ora{C} w$ 
is a cycle in $G$ such that $V(D) = V(C) \cup V(P)$, a contradiction. 
Thus $s$ is odd. 
Let 
$D = wx' \ora{P} xu u_{1} u^{+2} \ora{C} v_{2} u_{2} u_{3} v_{2}^{+} \ora{C} v_{4} u_{4}u_{5} v_{4}^{+} \ora{C} \dots v_{s-1}$ $u_{s-1}u_{s}v_{s-1}^{+} \ora{C} w$.  
Then $D$ is a cycle in $G$ such that $V(D) = (V(C) \bss \{u^{+}\}) \cup V(P)$. 
Hence by the maximality of $|V(C)|$, $V(P) = \{x\}$, in particular, $w \neq v$. 
Moreover, if there exists a vertex $y$ in $N_{G}(u^{+}; G-V(C))$, 
then by Claim~\ref{yu^{+3}}, $yu^{+3} \in E(G)$, 
and hence 
$(D - \{u^{+2}u^{+3} \}) + \{ u^{+}u^{+2}, yu^{+}, yu^{+3} \}$ is a longer cycle than $C$ 
(note that by Lemma~\ref{maximality of C} (i), $y \notin V(H)$), a contradiction.
Thus $N_{G}(u^{+}; G-V(C)) = \emptyset$. 
Therefore, $u^{+}$ is a component of $G-V(D)$, 
and $G-V(D)$ contains some components whose union is $H - \{x\}$. 
This implies that 
either 
$\mu(D) < \mu(C)$, or 
$\mu(D) = \mu(C)$ and 
$\omega(D) < \omega(C)$ holds, 
which contradicts 
the choice (C1) or (C2). 

This completes the proof of Theorem~\ref{tri th}.
\qed



\end{document}